# DISTRIBUTIONS OF LINEAR FUNCTIONALS OF TWO PARAMETER POISSON–DIRICHLET RANDOM MEASURES

By Lancelot F. James, Antonio Lijoi[1] and Igor Prünster[2]

*Hong Kong University of Science and Technology,
Università degli Studi di Pavia and
Università degli Studi di Torino*

The present paper provides exact expressions for the probability distributions of linear functionals of the two-parameter Poisson–Dirichlet process $\mathrm{PD}(\alpha,\theta)$. We obtain distributional results yielding exact forms for density functions of these functionals. Moreover, several interesting integral identities are obtained by exploiting a correspondence between the mean of a Poisson–Dirichlet process and the mean of a suitable Dirichlet process. Finally, some distributional characterizations in terms of mixture representations are proved. The usefulness of the results contained in the paper is demonstrated by means of some illustrative examples. Indeed, our formulae are relevant to occupation time phenomena connected with Brownian motion and more general Bessel processes, as well as to models arising in Bayesian nonparametric statistics.

**1. Introduction.** Let $(P_i)_{i\geq 1}$, with $P_1 > P_2 > \cdots > 0$ and $\sum_{k=1}^\infty P_k = 1$, denote a sequence of (random) ranked probabilities having the two-parameter $(\alpha,\theta)$ Poisson–Dirichlet law, denoted by $\mathrm{PD}(\alpha,\theta)$ for $0 \leq \alpha < 1$ and $\theta \geq 0$. A description, as well as a thorough investigation of its properties, is provided in [37]; see also [31, 32] and [35]. Equivalently, letting $V_k$, for any $k \geq 1$, denote independent random variables such that $V_k$ has $\mathrm{Beta}(1-\alpha,\theta+k\alpha)$ distribution, the $\mathrm{PD}(\alpha,\theta)$ law is defined as the ranked values of the stick-breaking sequence $W_1 = V_1$, $W_k = V_k \prod_{j=1}^{k-1}(1-V_j)$ for $k \geq 2$. Interestingly,

Received February 2007; revised June 2007.
[1]Supported by MUR Grant 2006/134525. Also affiliated with CNR-IMATI, Milan, Italy.
[2]Supported by MUR Grant 2006/133449. Also affiliated with Collegio Carlo Alberto and ICER, Turin, Italy.

PD$(\alpha,\theta)$ laws can also be obtained by manipulating random probabilities of the type $P_i = J_i/\tilde{T}$, where $\tilde{T} = \sum_{i=1}^{\infty} J_i$ and the sequence $(J_i)_{i\geq 1}$ stands for the ranked jumps of a subordinator. If the $J_i$'s are the ranked jumps of a gamma subordinator considered up to time $\theta > 0$, then the total mass $\tilde{T}$ has a gamma distribution with shape $\theta$ and scale 1, and $(P_i)_{i\geq 1}$ follows a PD$(0,\theta)$ law. At the other extreme, letting the $J_i$'s be the ranked jumps of a stable subordinator of index $0 < \alpha < 1$, $(P_i)_{i\geq 1}$ follows a PD$(\alpha,0)$ distribution. For both $\alpha$ and $\theta$ positive, the PD$(\alpha,\theta)$ model arises by first taking the ranked jumps governed by the stable subordinator conditioned on their total mass $\tilde{T}$ and then mixing over a power tempered stable law proportional to $t^{-\theta}f_\alpha(t)$, where $f_\alpha(t)$ denotes a density, with respect to the Lebesgue measure on $\mathbb{R}$, of an $\alpha$-stable random variable. We further recall that there is also the case of PD$(-\kappa, m\kappa)$, where $\kappa > 0$ and $m = 1, 2, \ldots$, which corresponds to symmetric Dirichlet random vectors of dimension $m$ and parameter $\kappa$. All of these models represent natural extensions of the important one-parameter family of Poisson–Dirichlet distributions, PD$(0,\theta)$, which is closely connected with the Dirichlet process.

Specifically, the corresponding PD$(\alpha,\theta)$ random probability measures are defined as follows. Independent of the sequence $(P_i)_{i\geq 1}$, or equivalently to the stick-breaking weights $(V_i)_{i\geq 1}$, let $(Z_i)_{i\geq 1}$ denote a collection of independent and identically distributed (i.i.d.) random elements which take values in a Polish space $\mathbb{X}$ endowed with the Borel $\sigma$-algebra $\mathscr{X}$ and have common nonatomic distribution $H$. One can then construct a PD$(\alpha,\theta)$ class of random probability measures, as

$$\tilde{P}_{\alpha,\theta}(\cdot) = \sum_{k=1}^{\infty} P_k \delta_{Z_k}(\cdot) = \sum_{k=1}^{\infty} W_k \delta_{Z_k}(\cdot).$$

When $\alpha = 0$, this is equivalent to the Dirichlet process which represents a cornerstone in Bayesian nonparametric statistics; see [10, 11, 12]. The law of $\tilde{P}_{\alpha,\theta}$ may be denoted as $\mathscr{P}_{(\alpha,\theta)}(\cdot|H)$. In particular, a random probability measure with distribution $\mathscr{P}_{(-\kappa,m\kappa)}(\cdot|H)$ can be represented as

$$\tilde{P}_{-\kappa,m\kappa}(\cdot) = \sum_{i=1}^{m} \frac{G_i}{\tilde{G}} \delta_{Z_i}(\cdot), \tag{1}$$

where $\tilde{G} = \sum_{i=1}^{m} G_i$ and the $G_i$'s are independent with Gamma$(\kappa,1)$ distribution, which in our notation, means that a density function for $G_i$ is of the form $[\Gamma(\kappa)]^{-1} x^{\kappa-1} e^{-x}$ for any $x > 0$. In [33], one can find a description of this model as Fisher's model; see also [19] for more references.

The choice of $\tilde{P}_{\alpha,\theta}$ for $\alpha > 0$, or of $\tilde{P}_{-\kappa,m\kappa}$, has attractive features which make them viable models for Bayesian nonparametric analysis, as shown in [5, 33], and [18, 19]. However, for the case $\alpha > 0$, most investigations of PD$(\alpha,\theta)$ laws appear in applications related to excursion/occupation time phenomena, as outlined in [37, 38] and, more recently, to combinatorial/probabilistic aspects of coalescent and phylogenetic processes; see [35] and [4] for numerous references along this line of research.

This paper studies the laws of mean functionals of the PD$(\alpha,\theta)$ class. We also briefly address the case PD$(-\kappa, m\kappa)$, which, as we shall show, essentially follows from the case of the Dirichlet process. In particular, for any nonnegative-valued function $f$ such that $\tilde{P}_{\alpha,\theta}(f)$ is finite, we obtain explicit formulae for the density and the cumulative distribution function (c.d.f.) of linear functionals

$$(2) \quad \tilde{P}_{\alpha,\theta}(f) = \int_{\mathbb{X}} f(x)\tilde{P}_{\alpha,\theta}(dx) = \sum_{k=1}^{\infty} P_k f(Z_k) = \sum_{k=1}^{\infty} f(Z_k) V_k \prod_{j=1}^{k-1}(1-V_j).$$

As such, we extend analogous formulae for Dirichlet processes, corresponding to the case of $\alpha = 0$, given by [6]. We do this by first resorting to the Cauchy–Stieltjes transforms of order $\theta$ for $\tilde{P}_{\alpha,\theta}(f)$, as derived in [44, 45], and also to a transform of order $\theta + 1$ deduced from [21], where, in particular, $\theta = 0$ for $\tilde{P}_{\alpha,0}(f)$. We then apply an Abel-type inversion formula described in [42] and finally combine those results with mixture representations of $\tilde{P}_{\alpha,\theta}(f)$ laws. We note that the case of $\tilde{P}_{\alpha,0}(f)$ for general $f$ is the most tractable, yielding explicit and simple expressions for the densities and c.d.f. which are expressed in terms of Abel transforms of $H$. The fact that our results have a strong connection to Abel transforms should not be totally surprising in view of the work in [14] where the laws of integrals of Bessel local times are investigated.

Interest in the results we are going to display and prove might arise in various contexts. In the paper, we will focus on specific issues related to (i) properties of the paths of Bessel processes and (ii) Bayesian nonparametric statistical inference. As for the former, we obtain results for pairs of parameters of the type $(\alpha, 0)$ and $(\alpha, \alpha)$, thus achieving useful expressions for the distribution of the lengths of excursions of Bessel processes and Bessel bridges. For example, we recover the important special cases of PD$(1/2, 0)$ and PD$(1/2, 1/2)$ which correspond to lengths of excursions of Brownian motion and Brownian bridge, respectively. As for (ii) above, knowledge of the probability distribution of $\tilde{P}_{\alpha,\theta}(f)$ can be useful for prior specification in applications where one is interested in making inference on a mean of the Poisson–Dirichlet process. However, apart from these two areas of research witch we are going to describe in more detail in the sequel, it is worth mentioning other potential applications of our results. For example, in [9, 25] and [45], it has been shown that results on means of the Dirichlet process have implications and interpretations relevant to, for example, the Markov moment problem, continued fraction theory and exponential representations of analytic functions. Since the $\mathrm{PD}(0,\theta)$ model can be seen as the limiting case of the $\mathrm{PD}(\alpha,\theta)$ distribution, as $\alpha \to 0$, we expect that some aspects of our work may also be applicable to these areas. Two other important applications of the $\mathrm{PD}(0,\theta)$ process for which our results could be of some interest relate to random combinatorial structures (see [2] for an exhaustive account) and population genetics.

The outline of the paper is as follows. In Section 2, we describe two areas of investigation for which our results are relevant and state an interesting distributional identity connecting $\mathrm{PD}(\alpha,\theta)$ and $\mathrm{PD}(0,\theta)$ means. In Section 3, we show how to use an inversion formula for Cauchy–Stieltjes transforms in order to determine a density function of $\tilde{P}_{\alpha,\theta}(f)$, as $(\alpha,\theta)$ varies in $(0,1) \times [0,+\infty)$. In Section 4, these general results are applied in order to determine generalized arcsine laws corresponding to mean functionals of a $\mathrm{PD}(\alpha,0)$ process: we show how to recover a well-known result and provide a representation for a density of $\int x \tilde{P}_{\alpha,0}(dx)$ when $\mathbb{E}[\tilde{P}_{\alpha,\theta}]$ coincides with the uniform distribution on $(0,1)$. Section 5 provides exact forms for a density of $\tilde{P}_{\alpha,\theta}(f)$ for any choice of $f$ which makes the random mean finite almost surely. In Section 6, a few distributional identities are given which prove to be useful in order to evaluate the distributions of means of $\mathrm{PD}(\alpha, 1-\alpha)$ and of $\mathrm{PD}(\alpha,\alpha)$ processes. Finally, Section 7 describes an algorithm for exact simulation whose formulation is suggested by results illustrated in Sections 2 and 6. All proofs are deferred to the Appendix.

We conclude this introductory section by recalling a useful fact. Indeed, note that $\tilde{P}_{\alpha,\theta}(f) \stackrel{d}{=} \int x \tilde{P}^*_{\alpha,\theta}(dx)$, where both $\tilde{P}_{\alpha,\theta}$ and $\tilde{P}^*_{\alpha,\theta}$ are Poisson–Dirichlet processes with $\mathbb{E}[\tilde{P}_{\alpha,\theta}(\cdot)] = H(\cdot)$ and $\mathbb{E}[\tilde{P}^*_{\alpha,\theta}(\cdot)] = H \circ f^{-1}(\cdot) =: \eta(\cdot)$. This explains why, with no loss of generality, we will confine ourselves to considering simple random means of the type $M_{\alpha,\theta}(\eta) := \int x \tilde{P}_{\alpha,\theta}(dx)$; see [41] for this line of reasoning. Moreover, the assumption of diffuseness of $H$ has been made only for consistency with the typical definition of the Poisson–Dirichlet random probability measure. Obviously, $\eta$ might have atoms, as will be seen in most of our examples; in any case, our results are still valid. Finally, in the sequel, $C_\eta$ will denote the convex hull of the support of $\eta$, that is, $C_\eta := \mathrm{co}(\mathrm{supp}(\eta))$.

**2. Related areas of application.** As already highlighted in the previous section, the main results achieved in the present paper find immediate application in two seemingly unrelated areas of research: the theory of Bessel processes and Bayesian nonparametric statistical inference. Below, we provide a brief description of the connection.

2.1. *Occupation times for Bessel processes and models for phylogenetic trees.* For functionals $\tilde{P}_{\alpha,\theta}(f)$, the generality of the space $\mathbb{X}$ is important as it allows one to formally describe phenomena where, for instance, $\mathbb{X}$ denotes path spaces of stochastic processes. Surprisingly, for general $(\alpha, \theta)$, very little is known about the laws of the simple, but important, case of $\tilde{P}_{\alpha,\theta}(\mathbb{I}_C)$ where $\mathbb{I}_C$ is the indicator function of set $C \in \mathscr{X}$, satisfying $\mathbb{E}[\tilde{P}_{\alpha,\theta}(\mathbb{I}_C)] = H(\mathbb{I}_C) = p \in (0,1)$. Hence, $f(Z) = \mathbb{I}_C(Z)$ is a Bernoulli random variable with success probability $p$, otherwise denoted Bernoulli($p$). In what follows, we will also let $\tilde{P}_{\alpha,\theta}(C)$ stand for the random mean $\tilde{P}_{\alpha,\theta}(\mathbb{I}_C)$ and $H(C) = H(\mathbb{I}_C)$. Using the representation provided in (2), one obtains

$$\tilde{P}_{\alpha,\theta}(C) = \sum_{k=1}^{\infty} Y_k V_k \prod_{j=1}^{k-1}(1 - V_j),$$

where $(Y_k)$ are i.i.d. Bernoulli($p$). The simple case corresponds to $\tilde{P}_{0,\theta}(C)$, which, since the Dirichlet process arises as a normalized gamma process, is well known to be a Beta($\theta p, \theta(1-p)$) random variable. This is apparent from the fact that if we set $Z_a$ to be a gamma random variable with density function at $x$ given by $[\Gamma(a)]^{-1} x^{a-1} e^{-x} \mathbb{I}_{\mathbb{R}^+}(x)$, then $\tilde{P}_{0,\theta}(C) \stackrel{d}{=} Z_{\theta p}/(Z_{\theta p} + Z_{\theta \bar{p}})$, where $Z_{\theta p}$ and $Z_{\theta \bar{p}}$ are independent and $\bar{p} = 1 - p$. The other known case corresponds to $\tilde{P}_{\alpha,0}(C)$, which has the Cauchy–Stieltjes transform

$$(3) \qquad \mathbb{E}[(1 + z\tilde{P}_{\alpha,0}(C))^{-1}] = \frac{(1+z)^{\alpha-1} p + \bar{p}}{(1+z)^{\alpha} p + \bar{p}}.$$

Such a transform has been inverted in [27], yielding, as $\alpha$ varies in $(0,1)$, the densities

$$(4) \qquad q_{\alpha,0}(x) = \frac{p\bar{p}\sin(\alpha\pi) x^{\alpha-1}(1-x)^{\alpha-1} \mathbb{I}_{(0,1)}(x)}{\pi[\bar{p}^2 x^{2\alpha} + p^2(1-x)^{2\alpha} + 2p\bar{p} x^{\alpha}(1-x)^{\alpha} \cos(\alpha\pi)]},$$

otherwise known as *generalized arcsine laws*. It is worth noting that this, as discussed in [3, 36] and [38], also corresponds to the fraction of time spent positive by a skew Bessel process of dimension $2 - 2\alpha$. Following [38], let $Y = (Y_t, t \geq 0)$ denote a real-valued continuous process such that (i) the zero set $Z$ of $Y$ is the range of an $\alpha$-stable subordinator and (ii) given $|Y|$, the signs of excursions of $Y$ away from zero are chosen independently of each other to be positive with probability $p$ and negative with probability $\bar{p} = 1 - p$. Examples of this kind of process are: Brownian motion ($\alpha = p = 1/2$); skew Brownian motion ($\alpha = 1/2$ and $0 < p < 1$); symmetrized Bessel process of dimension $2 - 2\alpha$; skew Bessel process of dimension $2 - 2\alpha$. Then, for any random time $T$ which is a measurable function of $|Y|$,

$$(5) \qquad A_T = \int_0^T \mathbb{I}_{(0,+\infty)}(Y_s)\,ds$$

denotes the time spent positive by $Y$ up to time $T$. Furthermore, remarkably, $A_T/T \stackrel{d}{=} A_t/t \stackrel{d}{=} A_1 = A$ and $A \stackrel{d}{=} \tilde{P}_{\alpha,0}(C)$. We see that the case of $\alpha = 1/2$, in (4) is the density found by [24] for the fraction of time spent positive by a skew Brownian motion. Moreover, when $p = 1/2$, this coincides with Lévy's famous result yielding the arcsine law for Brownian motion. That is, when $p = 1/2$, the random probability $\tilde{P}_{1/2,0}(C)$ has a Beta$(1/2, 1/2)$ distribution; see [28].

In [38], it is also shown that the fraction of time spent positive by a skew Bessel bridge of dimension $2 - 2\alpha$ corresponds to the law of $\tilde{P}_{\alpha,\alpha}(C)$. This random variable also arises, among other places, in Corollary 33 of [34]. Another recent instance is that of $\tilde{P}_{\alpha,1-\alpha}(C)$, which equates with the limiting distribution of a phylogenetic tree model described in Proposition 20 of [16]. However, results for these models are only well known for $\alpha = 1/2$, which corresponds to skew Brownian bridges. In particular, setting $p = 1/2$ yields the Lévy result for Brownian bridge which implies that $\tilde{P}_{1/2,1/2}(C)$ is uniform on $[0, 1]$. A density for $\tilde{P}_{1/2,\theta}(C)$ and general $p$ has been obtained by several authors; see, for instance, equation (3.4) in [5]. The case of $(1/2, \theta)$ when $p = 1/2$ is then Beta$(\theta + 1/2, \theta + 1/2)$; see also equation (65) in [1] for the density of $\tilde{P}_{1/2,1/2}(C)$ for general $p$ and yet another application related to the law of $\tilde{P}_{\alpha,\alpha}(C)$.

While the cases of Bernoulli $Y_k$'s are indeed quite interesting, we do wish to reiterate that it is substantially more difficult to obtain results for the more general case where the $Y_k$'s have a general distribution $\eta$.

2.2. *Bayesian nonparametric statistics.* The topic of this paper can be naturally connected to a large body of literature in Bayesian nonparametrics which is aimed at investigating the probability distribution of functionals of random probability measures. Besides the pioneering work in [6], we mention: [7] and [30], where nonlinear functionals of the Dirichlet process are studied; [29, 40] and [17], which provide developments and refinements of the earlier results in [6]; [21] and [41], which yield distributional results for means of a class in random probability measures that generalize the Dirichlet process. The interest in random probability measures in Bayesian nonparametric statistics is motivated by the fact that they define priors on spaces of probability distributions.

Here, we briefly describe the contribution in [6] since it has inspired our own approach. A first result contained in that paper consists of an important formula for the generalized Cauchy–Stieltjes transform of order $\theta$ of the mean functional $\tilde{P}_{0,\theta}(f)$ of the Dirichlet process $\tilde{P}_{0,\theta}$ with parameter measure $\theta H$. Supposing that $f : \mathbb{X} \to \mathbb{R}$ is such that $\int_{\mathbb{X}} \log(1 + |f(x)|) H(dx) < \infty$, the Cifarelli and Regazzini [6] show that

$$(6) \qquad \mathbb{E}\left[\frac{1}{(z + \tilde{P}_{0,\theta}(f))^\theta}\right] = e^{-\theta \int_{\mathbb{X}} \log(z + f(y)) H(dy)} = e^{-\theta \int_{\mathbb{R}} \log(z+x) \eta(dx)}$$

for any $z \in \mathbb{C}$ such that $|\arg(z)| < \pi$ and $\eta = H \circ f^{-1}$. The expression in (6) establishes that the Cauchy–Stieltjes transform of order $\theta$ of $\tilde{P}_{0,\theta}(f)$ is equivalent to the Laplace transform of $G_\theta(f)$, where $\tilde{P}_{0,\theta}(f) = G_\theta(f)/G_\theta(1)$ and $G_\theta$ is a gamma process with shape $\theta H$. The importance of (6) in different contexts was recognized by [9, 25] and [45]. In this regard, it is called the Markov–Krein identity for means of Dirichlet processes. It is called the Cifarelli–Regazzini identity in [22]. With considerable effort, Cifarelli and Regazzini [6] then apply an inversion formula to (6) to obtain an expression for the distribution of $\int x \tilde{P}_{0,\theta}(dx)$ as follows. Let $q_{0,\theta}$ denote the density of $\int x \tilde{P}_{0,\theta}(dx)$, where $\tilde{P}_{0,\theta}$ is a Dirichlet process with parameter measure $\theta \eta$, set $\Psi(x) := \eta((0,x])$ for any $x > 0$, and let

$$(7) \qquad R(t) = \int_{\mathbb{R}^+ \setminus \{t\}} \log(|t-x|) \eta(dx).$$

Then, from [6] (see also [7]), one has, for $\theta = 1$,

$$(8) \qquad q_{0,1}(x) = \frac{1}{\pi} \sin(\pi \Psi(x)) e^{-R(x)}$$

and when $\theta > 1$,

$$(9) \qquad q_{0,\theta}(x) = (\theta - 1) \int_0^x (x-t)^{\theta-2} \frac{1}{\pi} \sin(\pi \theta \Psi(t)) e^{-\theta R(t)}\, dt.$$

Additionally, an expression for the c.d.f., which holds for $\theta \Psi$ not having jumps greater than or equal to 1, is given by [6] as

$$(10) \qquad \int_0^x (x-t)^{\theta-1} \frac{1}{\pi} \sin(\pi \theta \Psi(t)) e^{-\theta R(t)}\, dt.$$

In particular, 10 holds for all $\theta > 0$ if $\eta$ is nonatomic. We note that while there are various formulae to describe the densities of $\int x \tilde{P}_{0,\theta}(dx)$, descriptions for the range $0 < \theta < 1$ prove to be difficult; see, for example, [6, 40] and [29].

Here, we provide a new description for the density, which holds for all $\theta > 0$. This result will be obvious from our subsequent discussion concerning the inversion formula for the Cauchy–Stieltjes transform and otherwise follows immediately from (10).

PROPOSITION 2.1. *Assume that $\eta$ admits a density on $\mathbb{R}^+$ and suppose that $R$ defined in (7) is differentiable. Then, the density of the Dirichlet process mean functional $\int x \tilde{P}_{0,\theta}(dx)$ may be expressed, for all $\theta > 0$, as*

$$(11) \qquad q_{0,\theta}(y) = \frac{1}{\pi} \int_0^y (y-t)^{\theta-1} d_{\theta,\eta}(t)\, dt,$$

*where*

$$(12) \qquad d_{\theta,\eta}(t) = \frac{d}{dt} \sin(\pi \theta \Psi(t)) e^{-\theta R(t)}.$$

It is apparent that practical usage of these formulae require tractable forms for $R$ and its derivative, which are not always obvious.

We close the present section with an important distributional identity which connects $\tilde{P}_{\alpha,\theta}(f)$ with mean functionals of the Dirichlet process. For a measurable function $f:\mathbb{X}\to\mathbb{R}^+$ such that $H(f^\alpha)<\infty$, let $Q_{\alpha,0}$ denote the probability distribution of $\tilde{P}_{\alpha,0}(f)$. This means that $\int_\mathbb{X} \log(1+x)Q_{\alpha,0}(dx) = \mathbb{E}[\log(1+\tilde{P}_{\alpha,0}(f))]$. The last expression is finite, since $\int_\mathbb{X} f^\alpha(x)H(dx)<\infty$. Hence, $\int_{\mathbb{R}^+} x\tilde{P}_{0,\theta}(dx)<\infty$ almost surely. From Theorem 4 in [45], note that, for any $z$ such that $|z|<\pi$,

$$(13) \quad \exp\left\{-\theta\int_{\mathbb{R}^+}\log[z+x]Q_{\alpha,0}(dx)\right\} = \left\{\int_{\mathbb{R}^+}[z+x]^\alpha\eta(dx)\right\}^{-\theta/\alpha}.$$

From the Cifarelli–Regazzini identity (6), the left-hand side in (13) coincides with the generalized Stieltjes transform of order $\theta$ of the random Dirichlet mean $M_{0,\theta}(Q_{\alpha,0})$, whereas the right-hand side is the generalized Stieltjes transform of order $\theta$ of $\tilde{P}_{\alpha,\theta}(f)$. These arguments can be summarized in the following theorem, which states a distributional identity between the mean $\int_{\mathbb{R}^+} x\tilde{P}_{0,\theta}(dx)$ of a Dirichlet process and a suitable linear functional of a $\mathrm{PD}(\alpha,\theta)$ process.

THEOREM 2.1. *Let $f:\mathbb{X}\to\mathbb{R}^+$ be a measurable function such that $H(f^\alpha)<\infty$, where $\alpha\in(0,1)$. If $Q_{\alpha,0}$ stands for the probability distribution of $\tilde{P}_{\alpha,0}(f)$, with $\tilde{P}_{\alpha,0}$ such that $\mathbb{E}[\tilde{P}_{\alpha,0}(\cdot)]=H(\cdot)$, then*

$$(14) \qquad \tilde{P}_{\alpha,\theta}(f) \stackrel{d}{=} \int_{\mathbb{R}^+} x\tilde{P}_{0,\theta}(dx),$$

*where $\tilde{P}_{0,\theta}$ is a Dirichlet process with $\mathbb{E}[\tilde{P}_{0,\theta}(\cdot)]=Q_{\alpha,0}(\cdot)$ and the symbol $\stackrel{d}{=}$ is used to denote equality in distribution. Moreover, the probability distribution of $\tilde{P}_{\alpha,\theta}(f)$ is absolutely continuous with respect to the Lebesgue measure on $\mathbb{R}$.*

It is worth noting that an alternative proof of (14) can be given, based on a construction of Gnedin and Pitman [13].

The last part of the statement in Theorem 2.1—that is the absolute continuity of the probability distribution of $\tilde{P}_{\alpha,\theta}(f)$ with respect to the Lebesgue measure—can be deduced from the absolute continuity of the probability distribution of $\int_{\mathbb{R}^+} x\tilde{P}_{0,\theta}(dx)$ which is proved in Proposition 2 of [29]. Hence, the probability distribution $Q_{\alpha,\theta}$ of $\tilde{P}_{\alpha,\theta}(f)$ has a density with respect to the Lebesgue measure on $\mathbb{R}$, which we denote $q_{\alpha,\theta}$.

**3. The probability distribution of $\tilde{P}_{\alpha,\theta}(f)$.** The present section provides a general expression of the density function of the mean of a two-parameter Poisson–Dirichlet distribution by resorting to an inversion formula for the (generalized) Cauchy–Stieltjes transform of order $\theta > 0$. Before getting into the details, let us introduce some new notation that will be used henceforth. Let $H$ be some nonatomic distribution on $(\mathbb{X}, \mathcal{X})$ and let $f : \mathbb{X} \to \mathbb{R}^+$ be any function in the set

$$\mathscr{E}_\alpha(H) := \{f : \mathbb{X} \to \mathbb{R}^+ \text{ s.t. } f \text{ is measurable and } H(f^\alpha) < +\infty\}.$$

Moreover, as anticipated in the Introduction, $\tilde{P}_{\alpha,\theta}$ denotes a random probability measure with law $\mathscr{P}_{(\alpha,\theta)}(\cdot|H)$. We confine our attention to functions in $\mathscr{E}_\alpha(H)$ for two reasons. First, the integrability condition $H(f^\alpha) = \int_\mathbb{X} f^\alpha(x) H(dx) < +\infty$ is necessary and sufficient for the (almost sure) finiteness of $\tilde{P}_{\alpha,0}(f)$; see Proposition 1 in [41] for a proof of this result. Hence, one can use the absolute continuity of $\mathscr{P}_{(\alpha,\theta)}(\cdot|H)$ with respect to $\mathscr{P}_{(\alpha,0)}(\cdot|H)$ in order to obtain $\tilde{P}_{\alpha,\theta}(f) < \infty$ with probability 1. Second, we consider only nonnegative functions since the inversion formula we resort to applies to functionals $\tilde{P}_{\alpha,\theta}(f)$ for a measurable function $f : \mathbb{X} \to \mathbb{R}^+$.

Given a function $f$ in $\mathscr{E}_\alpha(H)$, the transform of order $\theta > 0$ of $\tilde{P}_{\alpha,\theta}(f)$ is, for any $z \in \mathbb{C}$ such that $|\arg(z)| < \pi$,

$$(15) \quad \mathcal{S}_\theta[z; \tilde{P}_{\alpha,\theta}(f)] = \mathbb{E}\left[\frac{1}{(z + \tilde{P}_{\alpha,\theta}(f))^\theta}\right] = \left\{\int_\mathbb{X} [z + f(x)]^\alpha H(dx)\right\}^{-\theta/\alpha}.$$

Such a representation is to be attributed to [26] and also appears, with different proofs, in [44] and [45]. This transform turns out to work well in the case where $\theta > 1$. Additionally, we will need the transform of order $\theta + 1$, that is,

$$(16) \qquad \mathcal{S}_{\theta+1}[z; \tilde{P}_{\alpha,\theta}(f)] = \frac{\int_\mathbb{X} [z + f(x)]^{\alpha-1} H(dx)}{\{\int_\mathbb{X} [z + f(x)]^\alpha H(dx)\}^{\theta/\alpha+1}}.$$

In particular, for $\theta = 0$, we have, importantly, the Cauchy–Stieltjes transform of order 1 of the PD$(\alpha, 0)$ mean functionals,

$$(17) \qquad \mathcal{S}_1[z; \tilde{P}_{\alpha,0}(f)] = \frac{\int_\mathbb{X} [z + f(x)]^{\alpha-1} H(dx)}{\int_\mathbb{X} [z + f(x)]^\alpha H(dx)}.$$

The transforms (16) and (17) can be obtained as special cases of Proposition 6.2 in [21] with $n = 1$. Moreover, for $\theta > 0$, (16) can be obtained by taking a derivative of (15). When inverting (15) or (16), one obtains the probability distribution of $\tilde{P}_{\alpha,\theta}(f)$. Since Theorem 2.1 implies that the probability distribution of $\tilde{P}_{\alpha,\theta}(f)$ is, for any $(\alpha, \theta) \in (0, 1) \times \mathbb{R}^+$, absolutely continuous with respect to the Lebesgue measure on $\mathbb{R}$, we obtain the density function, $q_{\alpha,\theta}$ of $\tilde{P}_{\alpha,\theta}(f)$.

The particular inversion formula we are going to use has recently been given in [42]; see also [43]. In [8], one can find a detailed account of references on inversion formulae for generalized Cauchy–Stieltjes transforms.

THEOREM 3.1. *Let $f$ be a function in $\mathscr{E}_\alpha(H)$ and define*

$$(18) \quad \Delta_{\alpha,\theta}(t) := \frac{1}{2\pi\mathrm{i}} \lim_{\varepsilon\downarrow 0} \{\mathcal{S}_\theta[-t - \mathrm{i}\varepsilon; \tilde{P}_{\alpha,\theta}(f)] - \mathcal{S}_\theta[-t + \mathrm{i}\varepsilon; \tilde{P}_{\alpha,\theta}(f)]\}.$$

*The density function $q_{\alpha,\theta}$ of $\tilde{P}_{\alpha,\theta}(f)$, evaluated at a point $y$ in $C_\eta = \mathrm{co}(\mathrm{supp}(\eta))$, then coincides with*

$$(19) \quad q_{\alpha,\theta}(y) = \int_0^y (y - t)^{\theta-1} \Delta'_{\alpha,\theta}(t)\,dt.$$

*When $\theta > 1$, the expression above can be rewritten as*

$$(20) \quad q_{\alpha,\theta}(y) = (\theta - 1) \int_0^y (y - t)^{\theta-2} \Delta_{\alpha,\theta}(t)\,dt.$$

It is worth noting that if $\theta = 1$, then $q_{\alpha,1} = \Delta_{\alpha,1}$, thus yielding the same result as in [46]. The case corresponding to $\theta < 1$ can also be dealt with by computing the transform $\mathcal{S}_{\theta+1}[z; \tilde{P}_{\alpha,\theta}(f)]$. One then obtains

$$(21) \quad q_{\alpha,\theta}(y) = \theta \int_0^y (y - t)^{\theta-1} \tilde{\Delta}_{\alpha,\theta+1}(t)\,dt,$$

where

$$(22) \quad \begin{aligned} \tilde{\Delta}_{\alpha,\theta+1}(t) := \frac{1}{2\pi\mathrm{i}} \lim_{\varepsilon\downarrow 0} \bigg\{ & \frac{\int_{\mathbb{X}}[-t - \mathrm{i}\varepsilon + f(x)]^{\alpha-1} H(dx)}{[\int_{\mathbb{X}}(-t - \mathrm{i}\varepsilon + f(x))^\alpha H(dx)]^{\theta/\alpha+1}} \\ & - \frac{\int_{\mathbb{X}}[-t + \mathrm{i}\varepsilon + f(x)]^{\alpha-1} H(dx)}{[\int_{\mathbb{X}}(-t + \mathrm{i}\varepsilon + f(x))^\alpha H(dx)]^{\theta/\alpha+1}} \bigg\}. \end{aligned}$$

Note that the formulas (19) and (21) lead to the almost everywhere equality

$$(23) \quad \Delta'_{\alpha,\theta} = \theta \tilde{\Delta}_{\alpha,\theta+1}$$

for $\theta > 0$. Finally, note that $\tilde{\Delta}_{\alpha,1}$ is, by Widder's inversion, the density of $\tilde{P}_{\alpha,0}(f)$. Hence, a first approach for the determination of the distribution of $\tilde{P}_{\alpha,\theta}(f)$ will aim at the determination of $\Delta_{\alpha,\theta}$ and $\tilde{\Delta}_{\alpha,\theta+1}$. This task will be completed in the following sections. Note that once we obtain an explicit form for $\Delta_{\alpha,\theta}$, the c.d.f. of $\tilde{P}_{\alpha,\theta}(f)$ is given by

$$(24) \quad Q_{\alpha,\theta}((-\infty, x]) = \int_0^x q_{\alpha,\theta}(y)\,dy = \int_0^x (x - y)^{\theta-1} \Delta_{\alpha,\theta}(y)\,dy$$

for all $\theta > 0$. This result follows by using the representation in (19) and applying integration by parts. As we shall see, this representation plays a key role in obtaining simpler expressions and various identities for $\Delta_{\alpha,\theta}$, hence simplifying the formulas for the densities.

Finally, we anticipate some notation that will be useful. First, consider

$$\mathscr{A}^+_{\eta,d}(t) = \int_t^\infty (x-t)^d \eta(dx) \quad \text{and} \quad \mathscr{A}_{\eta,d}(t) = \int_0^t (t-x)^d \eta(dx),$$

which represent generalized Abel transforms of the measure $\eta$. Now, define

$$\gamma_\alpha(t) = \cos(\alpha\pi)\mathscr{A}_{\eta,\alpha}(t) + \mathscr{A}^+_{\eta,\alpha}(t),$$

$$\zeta_\alpha(t) = \sin(\alpha\pi)\mathscr{A}_{\eta,\alpha}(t), \qquad \rho_{\alpha,\theta}(t) = \frac{\theta}{\alpha}\arctan\frac{\zeta_\alpha(t)}{\gamma_\alpha(t)} + \frac{\pi\theta}{\alpha}\mathbb{I}_{\Gamma_\alpha}(t),$$

where $\Gamma_\alpha := \{t \in \mathbb{R}^+ : \gamma_\alpha(t) < 0\}$. Clearly, when $\alpha \leq 1/2$, $\gamma_\alpha(t) > 0$ for all $t$.

**4. Generalized arcsine laws.** We first deal with linear functionals of the $\mathrm{PD}(\alpha,0)$ process. As such, we generalize the results of Lamperti [27] for the case of $\tilde{P}_{\alpha,0}(C)$. We also point out that Regazzini, Lijoi and Prünster [41] obtain an expression for the c.d.f. of $M_{\alpha,0}(\eta)$ by exploiting a suitable inversion formula for characteristic functions and, additionally, provides expressions for its posterior density. Here, the approach we exploit leads to explicit and quite tractable expressions for the density which is expressed in terms of Abel transforms of $\eta$. Moreover, we also derive new expressions for the c.d.f. which can indeed be seen as generalized arcsine laws.

THEOREM 4.1. *Let $\eta$ be a probability measure on $(\mathbb{X},\mathscr{X})$ with $\mathbb{X} \subset \mathbb{R}^+$ and set $\Theta_{\alpha,\eta} := \{y \in \mathbb{R}^+ : \int_\mathbb{X} |y-t|^{\alpha-1}\eta(dt) < +\infty\} \cap C_\eta$. If $\int t^\alpha \eta(dt) < +\infty$ and the Lebesgue measure of $\Theta^c_{\alpha,\eta}$ is zero, then a density function of the random variable $M_{\alpha,0}(\eta) = \int x\tilde{P}_{\alpha,0}(dx)$, denoted by $q_{\alpha,0}$, coincides with*

$$(25) \quad q_{\alpha,0}(y) = \frac{\sin(\alpha\pi)}{\pi} \frac{\mathscr{A}^+_{\eta,\alpha}(y)\mathscr{A}_{\eta,\alpha-1}(y) + \mathscr{A}^+_{\eta,\alpha-1}(y)\mathscr{A}_{\eta,\alpha}(y)}{[\mathscr{A}^+_{\eta,\alpha}(y)]^2 + 2\cos(\alpha\pi)\mathscr{A}^+_{\eta,\alpha}(y)\mathscr{A}_{\eta,\alpha}(y) + [\mathscr{A}_{\eta,\alpha}(y)]^2}$$

*for any $y \in \Theta_{\alpha,\eta}$.*

The proof is provided in the Appendix. The result for the form of the density is new. We are also able to obtain, in view of obvious difficulties with direct integration, a rather remarkable expression of the c.d.f., given in the next theorem.

THEOREM 4.2. *Let $\eta$ be a probability measure on $\mathbb{R}^+$ such that $\int x^\alpha \eta(dx)$ is finite and the Lebesgue measure of the set $\Theta^c_{\alpha,\eta}$ is zero. If $t \mapsto \Psi(t) =$* $\eta((-\infty, t])$ *is Lipschitz of order* 1 *at any* $y \in \Theta_{\alpha,\eta}$, *then the c.d.f. of* $M_{\alpha,0}(\eta)$ *is given by*

$$(26) \qquad Q_{\alpha,0}((-\infty, x]) = \frac{1}{\alpha\pi} \arctan\left(\frac{\zeta_\alpha(x)}{\gamma_\alpha(x)}\right)$$

*for any* $x$ *in* $C_\eta$ *and* $\alpha$ *in* $(0,1)$.

Applying Theorem 4.2 to the case $\alpha = 1/2$, we obtain the following result.

COROLLARY 4.1. *Consider the setting as in Theorems 4.1 and 4.2. The density of the random variable* $M_{1/2,0}(\eta)$ *is given by*

$$(27) \qquad q_{1/2,0}(y) = \frac{1}{\pi} \frac{\mathscr{A}^+_{\eta,1/2}(y)\mathscr{A}_{\eta,-1/2}(y) + \mathscr{A}^+_{\eta,-1/2}(y)\mathscr{A}_{\eta,1/2}(y)}{[\mathscr{A}^+_{\eta,1/2}(y)]^2 + [\mathscr{A}_{\eta,1/2}(y)]^2}$$

*for any* $y \in \Theta_{\alpha,\eta}$ *and its c.d.f. is given by the generalized arcsine distribution*

$$(28) \quad Q_{1/2,0}((-\infty, x]) = \frac{2}{\pi} \arcsin(\pi \Delta_{1/2,1/2}(x)[\gamma^2_{1/2}(x) + \zeta^2_{1/2}(x)]^{1/2}).$$

REMARK 4.1. As we see, the results for the PD$(\alpha, 0)$ are tractable and, quite remarkably, only require the calculation of the Abel transforms $\mathscr{A}_{\eta,\alpha}$ and $\mathscr{A}^+_{\eta,\alpha}$. In this regard, one can, in general, obtain explicit results much more easily than for the case of the Dirichlet process. It is worth pointing out once again that our expressions for the c.d.f. show that these models have indeed generalized arcsine laws. These expression are rather surprising as it is not obvious how to integrate with respect to the densities.

Below, we illustrate a couple of applications of Theorem 4.1. The first example recovers a well-known result given in [27], while the second provides an expression for the density $q_{\alpha,0}$ when the parameter measure $\eta$ coincides with the uniform distribution on the interval $[0,1]$.

EXAMPLE 4.1 (Lamperti's occupation time density). Here, as a quick check of our results, we first revisit Lamperti's model. That is to say, the distribution of $\tilde{P}_{\alpha,0}(C)$. This corresponds to $\eta$ being the distribution of a Bernoulli distribution with success probability $p = \mathbb{E}[\tilde{P}_{\alpha,0}(C)] = 1 - \bar{p}$. It follows that for any $d > 0$, the Abel transforms for a Bernoulli random variable are given by

$$\mathscr{A}^+_{\eta,d}(t) = (1-t)^d p \quad \text{and} \quad \mathscr{A}_{\eta,d}(t) = t^d \bar{p}$$

for any $t$ in $(0,1)$. Hence, one easily sees that Lamperti's formula in (4) is recovered. In addition, we obtain the following new formula for the c.d.f.:

$$Q_{\alpha,0}((-\infty, x]) = \frac{1}{\alpha\pi} \arctan\left(\frac{\sin(\pi\alpha)x^\alpha \bar{p}}{\cos(\alpha\pi)x^\alpha \bar{p} + (1-x)^\alpha p}\right)$$

for any $x \in (0,1)$. This may also be expressed in terms of the arcsine, using the fact that for any $t \in (0,1)$,

$$\Delta_{\alpha,\alpha}(t) = \frac{\sin(\alpha\pi)t^\alpha \bar{p}}{\pi[t^{2\alpha}\bar{p}^2 + 2\cos(\alpha\pi)t^\alpha(1-t)^\alpha \bar{p}p + (1-t)^{2\alpha}p^2]}. \tag{29}$$

EXAMPLE 4.2 (Uniform parameter measure). We again note that, while there are several techniques one could have used to derive expressions for the functional $\tilde{P}_{\alpha,\theta}(C)$, it is considerably more difficult to obtain results for a more general choice of $\tilde{P}_{\alpha,\theta}(f)$, with $f$ in $\mathscr{E}_\alpha(H)$. Here, we demonstrate how our results easily identify the density in the case where $\eta(dx) = \mathbb{I}_{(0,1)}(x)\,dx$. For $M_{\alpha,0}(\eta) = \int x \tilde{P}_{\alpha,0}(dx)$, direct calculation of the Abel transforms leads to the expression of its density as

$$q_{\alpha,0}(y) = \frac{(\alpha+1)\sin(\alpha\pi)y^\alpha(1-y)^\alpha}{\alpha\pi[y^{2\alpha+2} + (1-y)^{2\alpha+2} + 2\cos(\alpha\pi)y^{\alpha+1}(1-y)^{\alpha+1}]}\mathbb{I}_{(0,1)}(y).$$

Note that one easily finds $\gamma_\alpha(t) = (t^{\alpha+1}\cos(\alpha\pi) + (1-t)^{\alpha+1})/(\alpha+1)$ and $\zeta_\alpha(t) = t^{\alpha+1}\sin(\alpha\pi)/(\alpha+1)$, additionally providing an expression for the c.d.f. In the Dirichlet case, the distribution of $\int_{(0,1)} x\tilde{P}_{0,\theta}(dx)$ can be determined by means of results contained in [6] and it is explicitly displayed in [9]. Its density function on $(0,1)$ has the form

$$q_{0,\theta}(y) = \frac{e}{\pi}(1-y)^{-(1+y)}y^{-y}\sin(\pi y)\mathbb{I}_{(0,1)}(y).$$

**5. The probability distribution of $\mathbf{PD}(\alpha,\theta)$ random means.** In the present section, we are going to consider more general cases than the $\alpha$-stable process dealt with in the previous section. In particular, we illustrate the results one can achieve by applying an inversion formula provided in [42] in order to obtain the density of $M_{\alpha,\theta}(\eta)$ when $\alpha \in (0,1)$ and $\theta > 0$ and then make a few remarks concerning the case in which $\alpha < 0$.

5.1. *Densities for $M_{\alpha,\theta}(\eta)$.* As suggested by the inversion formula, the evaluation of a density function for the random functional $M_{\alpha,\theta}$ basically amounts to the determination of the quantities $\Delta_{\alpha,\theta}$ and $\tilde{\Delta}_{\alpha,\theta+1}$. We then move on, stating the two main results of the section and providing an interesting illustration.

THEOREM 5.1. *For any $t \in \Theta_{\alpha,\eta}$ and $(\alpha,\theta) \in (0,1) \times \mathbb{R}^+$, one has*

$$\Delta_{\alpha,\theta}(t) = \frac{1}{\pi[\zeta_\alpha^2(t) + \gamma_\alpha^2(t)]^{\theta/(2\alpha)}}\sin(\rho_{\alpha,\theta}(t)). \tag{30}$$

A combination of the above result with Theorem 4.2 leads to another representation for the c.d.f. of the mean $\tilde{P}_{\alpha,0}(f)$ of an $\alpha$-stable subordinator. First, note that, for any $\alpha \in (0,1)$, one has

$$\rho_{\alpha,\alpha}(t) = \arctan\left(\frac{\zeta_\alpha(t)}{\gamma_\alpha(t)}\right) + \pi \mathbb{I}_{\Gamma_\alpha}(t) = \arcsin\left(\frac{\zeta_\alpha(t)\text{sign}(\gamma_\alpha(t))}{\sqrt{\zeta_\alpha^2(t) + \gamma_\alpha^2(t)}}\right) + \pi \mathbb{I}_{\Gamma_\alpha}(t),$$

where $\text{sign}(a) = 1$ if $a > 0$ and $\text{sign}(a) = -1$ if $a < 0$. At this point, from (30), one has

$$\Delta_{\alpha,\alpha}(t) = \frac{\text{sign}(\gamma_\alpha(t))}{\pi\sqrt{\zeta_\alpha^2(t) + \gamma_\alpha^2(t)}} \frac{\zeta_\alpha(t)|\gamma_\alpha(t)|}{\gamma_\alpha(t)\sqrt{\zeta_\alpha^2(t) + \gamma_\alpha^2(t)}} = \frac{\zeta_\alpha(t)}{\pi[\zeta_\alpha^2(t) + \gamma_\alpha^2(t)]}.$$

If, as in the statement of Theorem 4.2, one has that $\Psi$ is Lipschitz of order 1 at any $y \in \Theta_{\alpha,\eta}$, an alternative representation for the distribution function of a normalized $\alpha$-stable random mean displayed in (26), that is,

$$(31) \qquad Q_{\alpha,0}((-\infty,x]) = \frac{1}{\alpha\pi} \arcsin(\pi \Delta_{\alpha,\alpha}(x)\sqrt{\zeta_\alpha^2(x) + \gamma_\alpha^2(x)}),$$

holds true.

REMARK 5.1. As pointed out at the end of Section 3, if $\alpha \in (0,1/2]$, then $\rho_{\alpha,\theta}(t) = (\theta/\alpha) \arctan(\zeta_\alpha(t)/\gamma_\alpha(t))$. Hence, if one resorts to the expression for the c.d.f. of $\tilde{P}_{\alpha,0}(f)$ as provided in (31) it can be noted that $\rho_{\alpha,\theta}(t) = \theta\pi Q_{\alpha,0}((-\infty,t])$ and this yields a representation equivalent to the one displayed in (30), that is,

$$(32) \qquad \Delta_{\alpha,\theta}(t) = \frac{\sin(\theta\pi Q_{\alpha,0}((-\infty,t]))}{\pi[\zeta_\alpha^2(t) + \gamma_\alpha^2(t)]^{\theta/(2\alpha)}}.$$

A further evaluation of $\Delta_{\alpha,\theta}$ can be deduced by resorting to the correspondence between $\text{PD}(\alpha,\theta)$ means and Dirichlet means, as stated in Theorem 2.1. Indeed, one can prove the following useful identities.

THEOREM 5.2. Let $R_\alpha(t) = \int_{\mathbb{R}^+\setminus\{t\}} \log|t-y| Q_{\alpha,0}(dy)$. Then, for all $\theta > 0$, the following results hold:

(i) $\Delta_{\alpha,\theta}(t) = \pi^{-1} \sin(\pi\theta Q_{\alpha,0}((-\infty,t]))e^{-\theta R_\alpha(t)}$;
(ii) for any $t \in \{x : Q_{\alpha,0}((-\infty,x]) > 0\}$ and $\alpha \in (0,1)$,

$$e^{-R_\alpha(t)} = \left[\frac{\Delta_{\alpha,\alpha}(t)\pi}{\sin(\pi\alpha Q_{\alpha,0}(t))}\right]^{1/\alpha} = [\zeta_\alpha^2(t) + \gamma_\alpha^2(t)]^{-1/2\alpha}.$$

The expression of $\Delta_{\alpha,\theta}$, as determined in Theorem 5.1, Remark 5.1 or Theorem 5.2(i), is useful when one aims to evaluate the density function of $\tilde{P}_{\alpha,\theta}(f)$ corresponding to the case $\theta > 1$. When $\theta < 1$, one must resort to the expression given in (21) and then the evaluation of $\tilde{\Delta}_{\alpha,\theta+1}$ is necessary. The following theorem deals with this issue.

THEOREM 5.3. *For any $(\alpha, \theta) \in (0,1) \times \mathbb{R}^+$ and $t \in \Theta_{\alpha,\eta}$, the following identity holds true:*

$$\tilde{\Delta}_{\alpha,\theta+1}(t) = \frac{\gamma_{\alpha-1}(t)\sin(\rho_{\alpha,\theta}(t)) - \zeta_{\alpha-1}(t)\cos(\rho_{\alpha,\theta}(t))}{\pi[\zeta_\alpha^2(t) + \gamma_\alpha^2(t)]^{(\theta+\alpha)/(2\alpha)}}. \tag{33}$$

We now proceed to provide a simple illustration of the above results by means of an example. More detailed discussion about the determination of the probability distribution of $M_{\alpha,\theta}(\eta)$ is developed in the following section.

EXAMPLE 5.1 [First expressions for $\tilde{P}_{\alpha,\theta}(C)$]. It is interesting to compare the general case of $\tilde{P}_{\alpha,\theta}(C)$ with that of Lamperti's result in Example 4.1. Here, using the specifications in that example, it follows that

$$\Delta_{\alpha,\theta}(t) = \frac{\sin(\frac{\theta}{\alpha}\arctan(\frac{\bar{p}\sin(\alpha\pi)t^\alpha}{\bar{p}\cos(\alpha\pi)t^\alpha+p(1-t)^\alpha}) + \frac{\theta}{\alpha}\pi\mathbb{I}_{\Gamma_\alpha}(t))}{\pi\{\bar{p}^2 t^{2\alpha} + p^2(1-t)^{2\alpha} + 2\bar{p}p\cos(\alpha\pi)t^\alpha(1-t)^\alpha\}^{\theta/(2\alpha)}}, \tag{34}$$

where $\Gamma_\alpha = \varnothing$ if $\alpha \in (0, 1/2]$, whereas $\Gamma_\alpha = (0, \frac{v_\alpha}{1+v_\alpha})$ with $v_\alpha = (-p/(\bar{p}\cos(\alpha\pi)))^{1/\alpha}$ if $\alpha \in (1/2, 1)$. From (34), one can recover the expression for the c.d.f. of $\tilde{P}_{\alpha,\theta}(C)$ by resorting to (24) if the parameter $\theta > 1$. When $\theta < 1$, expressions for $\tilde{\Delta}_{\alpha,\theta+1}$ can also be calculated explicitly, leading to formulae for the density. In general, it is evident that such results are not as amenable as the case $\tilde{P}_{\alpha,0}(C)$, although they still lead to interesting insights. We will see that a case-by-case analysis can lead to more explicit expressions. We also develop other techniques in the forthcoming sections.

5.2. *The probability distribution of $M_{-\kappa,m\kappa}(\eta)$.* In this section, we establish the law of the mean functional of a random probability measure with distribution $\mathscr{P}_{(-\kappa,m\kappa)}(\cdot|H)$ which, according to (1), is given by

$$M_{-\kappa,m\kappa}(\eta) = \sum_{i=1}^m f(Z_i)\frac{G_i}{\widetilde{G}} = \sum_{i=1}^m Y_i \frac{G_i}{\widetilde{G}},$$

where the $Y_i$'s are i.i.d. with common probability distribution $\eta$. One reason to study these functionals is that for the choice of $\kappa = \theta/m$, one has that $M_{-\theta/m,\theta}(\eta)$ converges in distribution to $M_{0,\theta}(\eta)$ as $m \to \infty$. This fact may be found in, for example, [20]. It is easy to see that, conditionally on $(Y_1,\ldots,Y_m)$, $M_{-\kappa,m\kappa}(\eta) \stackrel{d}{=} M_{0,m\kappa}(\eta_m)$, where

$$\eta_m(\cdot) = \frac{1}{m}\sum_{i=1}^m \delta_{Y_i}(\cdot) \tag{35}$$

is the empirical distribution. Thus, descriptions of the conditional distribution of $M_{-\kappa,m\kappa}(\eta)$, given $(Y_i)_{i\geq 1}$, follow from (8), (9) and (10) for appropriate ranges of the parameter $\theta = m\kappa$, with $\eta$ replaced by $\eta_m$. The

Cauchy–Stieltjes transform of $M_{-\kappa,m\kappa}(\eta)$ of order $m\kappa$ is

$$\mathbb{E}\left[\frac{1}{(z+M_{-\kappa,m\kappa}(\eta))^{m\kappa}}\right] = \left[\int_0^\infty (z+y)^{-\kappa}\eta(dy)\right]^m, \qquad |\arg(z)| < \pi.$$

Now, set

$$\omega_m(t) := e^{-R(t)} = \prod_{i \in A_{t,m}} |t-y_i|^{-1},$$

where $A_{t,m} = \{i : y_i \neq t\} \cap \{1,\ldots,m\}$ and suppose that $c_{m\kappa}(t) := \int_{(0,\infty)} |t-y|^{-m\kappa}\eta(dy) < \infty$ for almost every $t$ with respect to the Lebesgue measure so that one can define $p_{m\kappa}(t) = \int_0^t (t-y)^{-m\kappa}\eta(dy)/c_{m\kappa}(t)$ and

$$h_{m,m\kappa}(t) := [c_{m\kappa}(t)]^m \sum_{j=1}^m \sin(\pi j\kappa) \binom{m}{j} [p_{m\kappa}(t)]^j [1-p_{m\kappa}(t)]^{m-j}.$$

This leads to the following interesting result.

PROPOSITION 5.1. *If $\mathbb{E}_\eta[\cdot]$ denotes the expected value taken with respect to $\eta$, then:*

(i) $\mathbb{E}_\eta[\omega_m^{m\kappa}(t)\sin(\pi m\kappa \eta_m(t))] = h_{m,m\kappa}(t)$;
(ii) *when $m\kappa = 1$, the density of $M_{-1/m,1}(\eta)$ is given by $h_{m,1}(x)/\pi$;*
(iii) *for $\kappa = \theta/m < 1$, the c.d.f. of $M_{-\theta/m,\theta}(\eta)$ is*

$$Q_{-\theta/m,\theta}((-\infty,x]) = \frac{1}{\pi}\int_0^x (x-t)^{\theta-1} h_{m,\theta}(t)\,dt.$$

*Furthermore, this c.d.f. converges to (10), as $m \to \infty$, for all $\theta > 0$.*

**6. Distributional recursions.** In this section, we describe mixture representations which are deducible from the posterior distribution of $\mathrm{PD}(\alpha,\theta)$ laws and existing results for the Dirichlet process. These represent aids in obtaining tractable forms of the distributions of various models $M_{\alpha,\theta}(\eta)$. In particular, we will use this to obtain a nice solution for all $\mathrm{PD}(\alpha,1-\alpha)$ models. Before stating the result, let us mention in advance that $B_{a,b}$ stands for a beta-distributed random variable with parameters $a, b > 0$.

THEOREM 6.1. *Let the random variables $W$, $M_{\alpha,\theta+\alpha}(\eta)$ and $B_{\theta+\alpha,1-\alpha}$ be mutually independent and such that $W$ has distribution $\eta$.*

(i) *Then, for $0 \leq \alpha < 1$ and $\theta > -\alpha$,*

$$M_{\alpha,\theta}(\eta) \stackrel{d}{=} B_{\theta+\alpha,1-\alpha} M_{\alpha,\theta+\alpha}(\eta) + (1 - B_{\theta+\alpha,1-\alpha})W.$$

*Note that when $\theta > 0$ and $\alpha = 0$, this equates with the mixture representation for Dirichlet process mean functionals.*

(ii) *Moreover, for any $\theta > 0$ and $\alpha \in (0,1)$, one has*

$$M_{\alpha,\theta}(\eta) \stackrel{d}{=} B_{\theta,1} M_{\alpha,\theta}(\eta) + (1 - B_{\theta,1}) M_{\alpha,0}(\eta).$$

An immediate consequence of this result is that if one has a tractable description of the distribution of $M_{\alpha,\theta+\alpha}(\eta)$, then one can easily obtain a description of the distribution of $M_{\alpha,\theta}(\eta)$. Moreover, from Theorem 6(i), one can deduce representations of $\tilde{P}_{\alpha,\theta}(f)$ for negative values of $\theta$, once an expression of $M_{\alpha,\theta+\alpha}(\eta)$ is available.

REMARK 6.1. Recall that $\tilde{P}_{1/2,0}(C)$ for $p = 1/2$ has the arcsine distribution Beta$(1/2, 1/2)$. Applying the mixture representation in statement (ii) of Theorem 6.1, one can see, via properties of Beta random variables, that $\tilde{P}_{1/2,\theta}(C)$ is Beta$(\theta + 1/2, \theta + 1/2)$. This corresponds to a result of ([7]) for $M_{0,\theta}(\eta)$, where $\eta$ is the arcsine law, although a connection to occupation time formula was not made there.

6.1. *Results for* PD$(\alpha,1)$ *and* PD$(\alpha,1-\alpha)$. We are now in a position to discuss some of the easiest (and also more important) cases. First, we deal with a $PD(\alpha,1)$ mean functional, $M_{\alpha,1}(\eta)$, and use it in order to determine, via a mixture representation, the probability distribution of $M_{\alpha,1-\alpha}(\eta)$.

We have already mentioned that when $\theta = 1$, the inversion formula simplifies and the density function of $M_{1,\alpha}(\eta)$ reduces to $q_{\alpha,1} = \Delta_{\alpha,1}$, as given in (30). For the range $0 < \alpha \leq 1/2$, this further reduces to

$$(36) \quad q_{\alpha,1}(y) = \Delta_{\alpha,1}(y) = \frac{\sin(1/\alpha \arcsin(\pi \Delta_{\alpha,\alpha}(t)\sqrt{\zeta_\alpha^2(y) + \gamma_\alpha^2(y)}))}{\pi[\zeta_\alpha^2(y) + \gamma_\alpha^2(y)]^{1/(2\alpha)}}.$$

If $\alpha = 1/n$, with $n = 2, 3, \ldots$, one can then use the multiple angle formula

$$\sin(nx) = \sum_{k=0}^{n} \binom{n}{k} [\cos(x)]^k [\sin(x)]^{n-k} \sin\left(\frac{\pi}{2}[n-k]\right),$$

in order to simplify the expression of $q_{\alpha,1}$. These remarks can be summarized as follows.

THEOREM 6.2. *A density function of $M_{\alpha,1}(\eta)$, for all $0 < \alpha < 1$, coincides with*

$$q_{\alpha,1}(y) = \Delta_{\alpha,1}(y) = \frac{\sin(1/\alpha \arctan \zeta_\alpha(y)/\gamma_\alpha(y) + \pi/\alpha \mathbb{I}_{\Gamma_\alpha}(y))}{\pi[\zeta_\alpha^2(y) + \gamma_\alpha^2(y)]^{1/2\alpha}}$$

*for any $y \in \Theta_{\alpha,\eta}$. When $\alpha = 1/n$, with $n = 2, 3, \ldots$, then a density function for $M_{1/n,1}(\eta)$ coincides with*

$$q_{1/n,1}(y) = \pi^{n-1} \Delta_{1/n,1/n}^n(y) \sum_{k=0}^{n} \binom{n}{k} \left(\frac{\gamma_{1/n}(y)}{\zeta_{1/n}(y)}\right)^k \sin\left(\frac{\pi}{2}[n-k]\right).$$

Setting $n=2$ in the last formula yields a simple form for a density function of $M_{1/2,1}(\eta)$, expressible as

$$q_{1/2,1}(y) = \frac{2}{\pi} \frac{\gamma_{1/2}(y)\zeta_{1/2}(y)}{[\zeta_{1/2}^2(y) + \gamma_{1/2}^2(y)]^2} = \frac{2}{\pi} \frac{\mathscr{A}_{\eta,1/2}(y)\mathscr{A}_{\eta,1/2}^+(y)}{[\zeta_{1/2}^2(y) + \gamma_{1/2}^2(y)]^2}.$$

Now, the density of $\mathrm{PD}(\alpha, 1-\alpha)$ mean functionals can be deduced from $\mathrm{PD}(\alpha, 1)$ models via the mixture representation given in Theorem 6.1.

THEOREM 6.3. *A density function of the random mean $M_{\alpha,1-\alpha}(\eta)$ is obtained via the distributional identity*

$$M_{\alpha,1-\alpha}(\eta) \stackrel{d}{=} B_{1,1-\alpha}M_{\alpha,1}(\eta) + (1 - B_{1,1-\alpha})W,$$

*where $B_{1,1-\alpha}$, $M_{\alpha,1}(\eta)$ and $W$ are independent. Here, $W$ is a random variable with distribution $\eta$. In particular, the density of $M_{\alpha,1-\alpha}(\eta)$ takes the form*

$$(1-\alpha)\int_0^\infty \int_0^1 \Delta_{\alpha,1}\left(\frac{x-wb}{1-b}\right)(1-b)^{-1}b^{-\alpha}\,db\,\eta(dw),$$

*where $\Delta_{\alpha,1}(t) = 0$ for any $t < 0$.*

EXAMPLE 6.1 (A distribution relevant to phylogenetic models). Recall from the Introduction that the random variable $\tilde{P}_{\alpha,1-\alpha}(C)$, when $\mathbb{E}[\tilde{P}_{\alpha,1-\alpha}(C)] = 1/2$, is equivalent in distribution to the random variable appearing in [16]. It is known that when $\alpha = 1/2$, the distribution is uniform, according to the well-known Lévy result; see [28]. Here, we obtain a quite tractable representation of the laws for all values of $\alpha$ and with $\mathbb{E}[\tilde{P}_{\alpha,1-\alpha}(C)] = p = 1 - \bar{p}$, for any $p \in (0,1)$. To this end, we first obtain the distribution of $\tilde{P}_{\alpha,1}(C)$. This can be easily obtained by setting $\theta = 1$ in (34), which yields the density function

$$f_{\alpha,1,p}(x) = \frac{\sin(\frac{1}{\alpha}\arctan(\frac{\bar{p}\sin(\alpha\pi)t^\alpha}{\bar{p}\cos(\alpha\pi)t^\alpha + p(1-t)^\alpha}) + \frac{\pi}{\alpha}\mathbb{I}_{\Gamma_\alpha}(t))}{\pi\{\bar{p}^2 t^{2\alpha} + p^2(1-t)^{2\alpha} + 2\bar{p}p\cos(\alpha\pi)t^\alpha(1-t)^\alpha\}^{1/(2\alpha)}},$$

where $\Gamma_\alpha = \varnothing$ if $\alpha \in (0, 1/2]$, whereas $\Gamma_\alpha = (0, \frac{v_\alpha}{1+v_\alpha})$ with $v_\alpha = (-p/(\bar{p}\cos(\alpha\pi)))^{1/\alpha}$ when $\alpha \in (1/2, 1)$. Since a density function $f_{\alpha,1,p}$ of $\tilde{P}_{\alpha,1}(C)$ is available, one can evaluate the probability distribution of $\tilde{P}_{\alpha,1-\alpha}(C)$ via the mixture representation stated in Theorem 6.1. It suffices to set $\eta = b_p$, where $b_p(x) = p^x(1-p)^{1-x}\mathbb{I}_{\{0,1\}}(x)$ is the probability mass functions of a Bernoulli random variable with parameter $p$. Hence, one has $\tilde{P}_{\alpha,1}(C) = M_{\alpha,1}(b_p)$ and $1 - M_{\alpha,1}(b_p) = 1 - \tilde{P}_{\alpha,1}(C) = M_{\alpha,1}(b_{\bar{p}})$.

COROLLARY 6.1. *Let $W$ denote a Bernoulli random variable with parameter $p$ and let $W$ be independent of $B_{1,1-\alpha}$ and $M_{\alpha,1}(b_p)$. Then, conditionally on the event $W = 1$, one has $M_{\alpha,1-\alpha}(b_p) \stackrel{d}{=} 1 - B_{1,1-\alpha}M_{\alpha,1}(b_{\bar{p}})$. On the other hand, given $W = 0$, one has $M_{\alpha,1-\alpha}(b_p) \stackrel{d}{=} B_{1,1-\alpha}M_{\alpha,1}(b_p)$. Equivalently, a density function of $M_{\alpha,1-\alpha}(b_p)$ is obtained via the distributional relationship*

$$f_{\alpha,1-\alpha,p}(t) = (1-\alpha)\int_0^1 \left[\bar{p}f_{\alpha,1,p}\left(\frac{t}{u}\right) + pf_{\alpha,1,\bar{p}}\left(\frac{1-t}{u}\right)\right]u^{-1}(1-u)^{-\alpha}\,du.$$

6.2. *The case of* $\mathrm{PD}(\alpha,\alpha)$. The important case of $\mathrm{PD}(\alpha,\alpha)$ is, in general, more challenging than the case of $\mathrm{PD}(\alpha,1-\alpha)$. Of course, these two agree in the case of $\alpha = 1/2$ corresponding to quantities related to Brownian bridges. Technically, one can apply the formula based on $\tilde{\Delta}_{\alpha,\alpha+1}$, but this does not always yield very nice expressions. Alternatively, in the special case where $1 - \alpha = 2\alpha$, that is $\alpha = 1/3$, one might think of using mixture representation results such as those given in Theorems 6.1 and 6.2. According to the latter, one can determine $M_{\alpha,1-\alpha}(\eta)$ and then, by resorting to the former (with $\theta = \alpha$), one obtains $M_{\alpha,\alpha}(\eta)$. For example, one can use the Dirichlet process mixture representation to obtain the probability distribution of $M_{2/3,2/3}(\eta)$ from the distribution of $M_{\frac{2}{3},\frac{1}{3}}(\eta)$. Additionally, when $\alpha > 1/2$, one may use the density representation of $M_{\alpha,2\alpha}(\eta)$ based on $\Delta_{\alpha,2\alpha}$, coupled with the mixture representation. Let us investigate these cases by considering specific examples.

EXAMPLE 6.2 [Probability distribution of $\tilde{P}_{\alpha,\alpha}(C)$]. First, note that, having set $p = 1 - \bar{p} = \eta(C) \in (0,1)$, the quantity $\Delta_{\alpha,\alpha}$ has been described in (29). The expressions for $\tilde{\Delta}_{\alpha,\alpha+1}$ are the same for any $\alpha \in (0,1)$ since $\sin(2\arctan(\frac{\zeta_\alpha(t)}{\gamma_\alpha(t)})) = \sin(2\arctan(\frac{\zeta_\alpha(t)}{\gamma_\alpha(t)}) + 2\pi\mathbb{I}_{\Gamma_\alpha}(t))$. Hence, for any $\alpha \in (0,1)$ and $t \in (0,1)$, one has

$$(37) \quad \tilde{\Delta}_{\alpha,\alpha+1}(t) = \frac{2\gamma_\alpha(t)\gamma_{\alpha-1}(t)\zeta_\alpha(t) - \zeta_{\alpha-1}(t)\gamma_\alpha^2(t) + \zeta_{\alpha-1}(t)\zeta_\alpha^2(t)}{\pi[\gamma_\alpha^2(t) + \zeta_\alpha^2(t)]^2},$$

with $\gamma_\alpha(t) = p(1-t)^\alpha + \bar{p}\cos(\alpha\pi)t^\alpha$ and $\zeta_\alpha(t) = \bar{p}\sin(\alpha\pi)t^\alpha$. These findings, with some simple algebra, lead to the following corollary.

COROLLARY 6.2. *The random probability $\tilde{P}_{\alpha,\alpha}(C)$ admits a density function coinciding with*

$$q_{\alpha,\alpha}(y) = \frac{\alpha\bar{p}\sin(\alpha\pi)}{\pi}\int_0^y [t(y-t)]^{\alpha-1}$$
$$\times (p^2(1-t)^{2\alpha-1}(1+t)$$

$$\begin{aligned}(38)\qquad &+ 2p\bar{p}t^{\alpha+1}(1-t)^{\alpha-1}\cos(\alpha\pi) - \bar{p}^2 t^{2\alpha})\\ &\times ([p^2(1-t)^{2\alpha} + \bar{p}^2 t^{2\alpha}\\ &\quad + 2p\bar{p}t^{\alpha}(1-t)^{\alpha}\cos(\alpha\pi)]^2)^{-1}\,dt\end{aligned}$$

for any $y$ in $(0,1)$.

It is now worth noting that the above formula, with $\alpha = p = 1/2$, yields the well-known result about the probability distribution of $A = \int_0^1 \mathbb{I}_{(0,+\infty)}(Y_s)\,ds$, in the case where the Markov process $Y = \{Y_t : t \in [0,1]\}$ is a Brownian bridge. Indeed, Lévy found that $A$ is uniformly distributed on the interval $(0,1)$; see [28]. In this case, $\tilde{\Delta}_{1/2,3/2}(t) = 2\pi^{-1}t^{-1/2}$ and the density function of $\tilde{P}_{1/2,1/2}(C)$ is given by

$$q_{1/2,1/2}(y) = \frac{1}{2\pi} 2 \int_0^y t^{-1/2}(y-t)^{-1/2}\,dt = 1$$

for any $y \in (0,1)$.

EXAMPLE 6.3 (Uniform parameter measure). Let us again consider the case in which $\eta(dx) = \mathbb{I}_{(0,1)}(dx)$. Recall that $\gamma_\alpha(t) = (t^{\alpha+1}\cos(\alpha\pi) + (1-t)^{\alpha+1})/(\alpha+1)$ and $\zeta_\alpha(t) = t^{\alpha+1}\sin(\alpha\pi)/(\alpha+1)$. These yield

$$\begin{aligned}\tilde{\Delta}_{\alpha,\alpha+1}(t) &= (1+\alpha)^2 \frac{\sin(\alpha\pi)}{\alpha\pi}\\ &\quad \times \frac{t^\alpha[(1-t)^{2\alpha+1}(1+t) - t^{2\alpha+2} + 2\cos(\alpha\pi)t^{\alpha+2}(1-t)^\alpha]}{[t^{2\alpha+2} + (1-t)^{2\alpha+2} + 2\cos(\alpha\pi)t^{\alpha+1}(1-t)^{\alpha+1}]^2}.\end{aligned}$$

The expression of the density $q_{\alpha,\alpha}$ somewhat simplifies when $\alpha = 1/2$. Indeed, in this case, one has

$$\tilde{\Delta}_{1/2,3/2}(t) = \frac{9\sqrt{t}[(1-t)^2(1+t) - t^3]}{2\pi[1 - 3t(1-t)]^2}$$

for any $t \in (0,1)$. In order to determine the probability density function $q$, compute

$$\begin{aligned}I_{r,s}(y) &:= \frac{2}{\pi} \int_0^y \frac{(y-t)^{-1/2}t^{r+1/2}(1-t)^s}{[1-3t(1-t)]^2}\,dt\\ &= \frac{2}{\pi} \sum_{n\geq 0} \frac{(2)_n 3^n}{n!} \int_0^y (y-t)^{-1/2}t^{r+n+1/2}(1-t)^{n+s}\,dt\\ &= \frac{2}{\pi} \sum_{n\geq 0} \frac{(2)_n 3^n}{n!} \sum_{k=0}^{n+s} \binom{n+s}{k}(-1)^k \int_0^y (y-t)^{-1/2}t^{r+n+k+1/2}\,dt\end{aligned}$$

$$= \frac{2}{\sqrt{\pi}} \sum_{n \geq 0} \frac{(2)_n 3^n}{n!} \sum_{k=0}^{n+s} \binom{n+s}{k} (-1)^k y^{n+k+r+1} \frac{\Gamma(r+n+k+3/2)}{\Gamma(r+n+k+2)},$$

where $(a)_n = \Gamma(a+n)/\Gamma(a)$ for any $a > 0$ and $n \geq 0$. Consequently,

$$q_{1/2,1/2}(y) = \int_0^y (y-t)^{-1/2} \tilde{\Delta}_{1/2,3/2}(t)\, dt = I_{0,2}(y) + I_{1,2}(y) - I_{1,0}(y)$$

for any $y$ in $(0,1)$.

An alternative representation of this density can be achieved by resorting to Theorem 6.1. Indeed, one has that $M_{1/2,1/2}(\eta) \stackrel{d}{=} B_{1,1/2} M_{1/2,1}(\eta) + (1 - B_{1,1/2})Y$, where the density function of $M_{1/2,1}(\eta)$ is given by

$$q_{1/2,1}(y) = \Delta_{1/2,1}(y) = \frac{9}{2\pi} \frac{y^{3/2}(1-y)^{3/2}}{\{y^3 + (1-y)^3\}^2} \mathbb{I}_{(0,1)}(y)$$

and $Y$ is uniformly distributed over the interval $(0,1)$. This then suggests that a density of $M_{1/2,1/2}(\eta)$ can be represented as

$$\begin{aligned}
q_{1/2,1/2}(y) &= \frac{1}{2} \int_0^{x_1} (x_1 - x_3)^{-1/2} \left\{ \int_{x_1}^1 (x_2 - x_3)^{-1/2} q_{1/2,1}(x_2)\, dx_2 \right\} dx_3 \\
&\quad + \frac{1}{2} \int_{x_1}^1 (x_3 - x_1)^{-1/2} \left\{ \int_0^{x_1} (x_3 - x_2)^{-1/2} q_{1/2,1}(x_2)\, dx_2 \right\} dx_3 \\
&= \frac{9}{2\pi} \sqrt{x_1} \int_{x_1}^1 \frac{x_2(1-x_2)^{3/2}}{\{x_2^3 + (1-x_2)^3\}^2}\, {}_2F_1\left(\frac{1}{2}, 1; \frac{3}{2}; \frac{x_1}{x_2}\right) dx_2 \\
&\quad + \frac{9}{2\pi} \sqrt{1-x_1} \int_0^{x_1} \frac{x_2^{3/2}(1-x_2)}{\{x_2^3 + (1-x_2)^3\}^2}\, {}_2F_1\left(\frac{1}{2}, 1; \frac{3}{2}; \frac{1-x_1}{1-x_2}\right) dx_2,
\end{aligned}$$

where ${}_2F_1$ is the Gauss hypergeometric function.

**7. Perfect sampling $M_{\alpha,\theta}(\eta)$.** Our results thus far have provided quite a few expressions for the densities and c.d.f.'s of $M_{\alpha,\theta}(\eta)$ which are certainly interesting from an analytic viewpoint. However, it is clear that if one were interested in drawing random samples, it is not always obvious how to do so. The clear exception for all $\eta$ is the $M_{\alpha,0}(\eta)$ case, where one can apply straightforward rejection sampling based on the explicit density in Theorem 4.1. Here, we show that this fact, in conjunction with the correspondence to the Dirichlet process established in Theorems 2.1 and 6.1, allows us to perfectly sample random variables $M_{\alpha,\theta}(\eta)$ for all $0 < \alpha < 1$ and $\theta > 0$. One simply needs to apply the perfect sampling procedure for Dirichlet mean functionals devised by [15]; see also [23] for an application of this idea to a class of non-Gaussian Ornstein–Uhlenbeck models arising in financial econometrics. First, recall that Theorem 6.1 establishes the distributional identity

$$M_{\alpha,\theta}(\eta) \stackrel{d}{=} M_{0,\theta}(Q_{\alpha,0}) \stackrel{d}{=} M_{0,\theta}(Q_{\alpha,0})B_{\theta,1} + (1-B_{\theta,1})M_{\alpha,0}(\eta).$$

Recognizing this, we now recount the basic elements of the perfect sampling algorithm of [15], tailored to the present situation. First, note that perfect sampling can be achieved if $0 \leq a \leq M_{\alpha,\theta}(\eta) \leq b < \infty$ almost surely. Furthermore, note that this is true if and only if the support of $Q_{\alpha,0}$ is $[a,b]$ or, equivalently, $M_{\alpha,0}(\eta) \in [a,b]$. Now, as explained in [15], following the procedure of [39], one can design an upper and lower dominating chain starting at some time $-N$ in the past up to time 0. The upper chain, say $uM_{\alpha,\theta}(\eta)$, is started at $uM_{\alpha,\theta,-N}(\eta) = b$ and the lower chain, $lM_{\alpha,\theta}(\eta)$, is started at $lM_{\alpha,\theta,-N}(\eta) = a$. One runs the Markov chains for each $n$ based on the equations

$$uM_{\alpha,\theta,n+1}(\eta) = B_{n,\theta}X_n + (1-B_{n,\theta})uM_{\alpha,\theta,n}(\eta)$$

and

$$lM_{\alpha,\theta,n+1}(\eta) = B_{n,\theta}X_n + (1-B_{n,\theta})lM_{\alpha,\theta,n}(\eta),$$

where the chains are coupled using the same random independent pairs $(B_{n,\theta}, X_n)$, where for each $n$, $B_{n,\theta}$ has a Beta$(1,\theta)$ distribution and $X_n$ has distribution $F_{\alpha,\eta}$. That is, $X_n \stackrel{d}{=} M_{\alpha,0}(\eta)$. The chains are said to *coalesce* when $D = |uM_{\alpha,\theta,n}(\eta) - lM_{\alpha,\theta,n}(\eta)| < \varepsilon$ for some small $\varepsilon$. Notice, importantly, that this method only requires the ability to sample $M_{\alpha,0}(\eta)$, which is provided by Theorem 4.1, and an independent Beta random variable.

## APPENDIX A: PROOF OF THEOREM 3.1

**A.1. An inversion formula for the Cauchy–Stieltjes transform.** In order to determine the density function, $q_{\alpha,\theta}$, of the random mean $\tilde{P}_{\alpha,\theta}(f)$, one can invert the transform $\mathcal{S}_\theta[z; \tilde{P}_{\alpha,\theta}(f)]$. Indeed, since $q_{\alpha,\theta}$ is a density function, with respect to the Lebesgue measure, on $\mathbb{R}^+$, one has $\int_0^c y^{-\rho} q_{\alpha,\theta}(y)\,dy < \infty$ for some $\rho > 0$ for any $c > 0$. The inversion formula provided in [42] can be applied to obtain

$$(39) \qquad q_{\alpha,\theta}(y) = -\frac{y^\theta}{2\pi \mathrm{i}} \int_{\mathcal{W}} (1+w)^{\theta-1} \mathcal{S}'_\theta[yw; \tilde{P}_{\alpha,\theta}(f)]\,dw.$$

In the previous formula, $\mathcal{W}$ is a contour in the complex plane starting at $w = -1$, encircling (in the counterclockwise sense) the origin and ending at $w = -1$, while $\mathcal{S}'_\theta[yw; \tilde{P}_{\alpha,\theta}(f)] = \frac{d}{dz}\mathcal{S}_\theta[z; \tilde{P}_{\alpha,\theta}(f)]|_{z=yw}$. If $\theta > 1$, one can integrate by parts, thus obtaining

$$(40) \qquad q_{\alpha,\theta}(y) = \frac{\theta-1}{2\pi \mathrm{i}} y^{\theta-1} \int_{\mathcal{W}} (1+w)^{\theta-2} \mathcal{S}_\theta[yw; \tilde{P}_{\alpha,\theta}(f)]\,dw;$$

see (18) and (19) in [42].

**A.2. Proof of Theorem 3.1.** In order to obtain (19) from the above representation (39), first note that the complex integral can be rewritten as follows. Set, for any $y$ in the convex hull of the support of $H \circ f^{-1}$, $\Gamma_y \subset \mathbb{C}$ to be the path starting at $w = -y$, encircling the origin and ending at $w = -y$. Hence, a change of variable in (39) yields

$$q_{\alpha,\theta}(y) = -\frac{1}{2\pi\mathrm{i}} \int_{\Gamma_y} (y+w)^{\theta-1} \mathcal{S}'_\theta[w; \tilde{P}_{\alpha,\theta}(f)] \, dw.$$

For simplicity, introduce the function $w \mapsto h(w) = (y+w)^{\theta-1} \mathcal{S}'_\theta[w; \tilde{P}_{\alpha,\theta}(f)]$ and note that $\mathcal{S}'_\theta[w; \tilde{P}_{\alpha,\theta}(f)] = \mathcal{S}_{\theta+1}[w; \tilde{P}_{\alpha,\theta}(f)]$, as described in (16). By virtue of Cauchy's theorem, one has

$$\int_{\Gamma_y} h(w) \, dw = \int_{-y+\mathrm{i}0}^{-y-\mathrm{i}\varepsilon} h(w) \, dw + \int_{-y}^{0} h(x-\mathrm{i}\varepsilon) \, dx$$
$$+ \mathrm{i}\varepsilon \int_{-\pi/2}^{\pi/2} e^{\mathrm{i}s} h(\varepsilon e^{\mathrm{i}s}) \, ds + \int_{0}^{-y} h(x+\mathrm{i}\varepsilon) \, dx + \int_{-y+\mathrm{i}\varepsilon}^{-y+\mathrm{i}0} h(w) \, dw.$$

A few straightforward simplifications lead to

(41)
$$\int_{\Gamma_y} h(w) \, dw = \int_0^y [h(-x-\mathrm{i}\varepsilon) - h(-x+\mathrm{i}\varepsilon)] \, dx$$
$$+ \mathrm{i}\varepsilon \int_{-\pi/2}^{\pi/2} e^{\mathrm{i}s} h(\varepsilon e^{\mathrm{i}s}) \, ds + \int_{-y+\mathrm{i}0}^{-y-\mathrm{i}\varepsilon} h(w) \, dw$$
$$+ \int_{-y+\mathrm{i}\varepsilon}^{-y+\mathrm{i}0} h(w) \, dw.$$

In order to determine the behavior, as $\varepsilon \downarrow 0$, of the last two summands in (41), let us show that the function $s \mapsto h(-y+\mathrm{i}s)$ is integrable in a neighborhood of the origin. Indeed, one has

$$|h(-y+\mathrm{i}s)| = |s|^{\theta-1} \frac{|\int_{\mathbb{X}} (-y+\mathrm{i}s + f(x))^{\alpha-1} H(dx)|}{|\{\int_{\mathbb{X}} (-y+\mathrm{i}s + f(x))^\alpha H(dx)\}^{\theta/\alpha+1}|}.$$

As for the numerator, one has $|\int_{\mathbb{X}} (-y+\mathrm{i}s + f(x))^{\alpha-1} H(dx)| \leq \int_{\mathbb{R}^+} |-y+x+\mathrm{i}s|^{\alpha-1} \eta(dx) \leq \int_{\mathbb{R}^+} |x-y|^{\alpha-1} \eta(dx)$ and this is finite for any $y$ in $\Theta_{\alpha,\eta}$. On the other hand, if one sets $g_1(x,y;s) := |x-y+\mathrm{i}s|^\alpha \cos(\alpha \arg(x-y+\mathrm{i}s))$ and $g_2(x,y;s) := |x-y+\mathrm{i}s|^\alpha \sin(\alpha \arg(x-y+\mathrm{i}s))$, the denominator would coincide with

$$\left| \left\{ \int_{\mathbb{R}^+} (x-y-\mathrm{i}s)^\alpha \eta(dx) \right\}^{\theta/\alpha+1} \right|$$
$$= \left| \exp\left\{ (\theta\alpha^{-1}+1) \log\left( \int_{\mathbb{R}^+} [g_1(x,y;s) + \mathrm{i}g_2(x,y;s)] \eta(dx) \right) \right\} \right|$$

$$= \left\{ \left( \int_{\mathbb{R}^+} g_1(x,y;s)\eta(dx) \right)^2 + \left( \int_{\mathbb{R}^+} g_2(x,y;s)\eta(dx) \right)^2 \right\}^{\theta(2\alpha)^{-1}+2^{-1}}$$

$$\geq \left( \int_{\mathbb{R}^+} g_2(x,y;s)\eta(dx) \right)^{\theta\alpha^{-1}+1}.$$

We now confine ourselves to the case where $s$ is in $(0,\xi)$ for some positive $\xi$. The same reasoning can be applied in the opposite case when $s \in (-\xi, 0)$. Since $\arg(x - y + \mathrm{i}s) = \arctan(s/(x-y)) + \pi \mathbb{I}_{(0,y)}(x)$, this quantity is a value in $(\pi/2, \pi)$ and, hence, $\sin(\alpha \arg(x - y + \mathrm{i}s)) > 0$. Consequently, if we set $w_\alpha = \alpha[\pi \mathbb{I}_{(0,2/3]}(\alpha) + (\pi/2)\mathbb{I}_{(2/3,1)}(\alpha)]$, the denominator is bounded below by

$$\left\{ \int_{\mathbb{R}^+} |x-y|^\alpha \sin\left( \alpha \arctan \frac{s}{x-y} + \alpha\pi \mathbb{I}_{(0,y)}(x) \right) \eta(dx) \right\}^{\theta\alpha^{-1}+1}$$

$$\geq \left\{ \sin(w_\alpha) \int_{(0,y]} |x-y|^\alpha \eta(dx) \right\}^{\theta\alpha^{-1}+1}$$

and the latter must be greater than some positive constant, say $K_y$, for any $y$ in $\Theta_{\alpha,\eta}$. Finally, we have $|h(-y + \mathrm{i}s)| \leq K'_y s^{\theta-1}$ for any $s$ in $(0,\xi)$. This implies that, when we let $\varepsilon \downarrow 0$, the last two terms in the right-hand side of (41) tend to zero. As far as the second summand in (41) is concerned, one can also show that it tends to zero since

$$(42) \qquad \lim_{\varepsilon \downarrow 0} \varepsilon^\alpha \int_{-\pi/2}^{\pi/2} e^{\mathrm{i}s}(y + \varepsilon e^{\mathrm{i}s})^{\theta-1} \varepsilon^{1-\alpha} \mathcal{S}_{\theta+1}[\varepsilon e^{\mathrm{i}s}; \tilde{P}_{\alpha,\theta}(f)]\, ds = 0.$$

In order to show this, observe that, for any real number $s$ such that $|s| \leq \pi/2$, one has $|(y + \varepsilon e^{-\mathrm{i}s})|^{\theta-1} = \{y^2 + \varepsilon^2 + 2y\varepsilon\cos(s)\}^{(\theta-1)/2} \leq y^{\theta-1}$, when $\theta < 1$. Moreover, one can determine a bound for

$$(43) \qquad \varepsilon^{1-\alpha} |\mathcal{S}_{\theta+1}[\varepsilon e^{\mathrm{i}s}; \tilde{P}_{\alpha,\theta}(f)]| = \varepsilon^{1-\alpha} \frac{|\int_{\mathbb{R}^+} [\varepsilon e^{\mathrm{i}s} + x]^{\alpha-1} \eta(dx)|}{|\{\int_{\mathbb{R}^+} [\varepsilon e^{\mathrm{i}s} + x]^\alpha \eta(dx)\}^{\theta/\alpha+1}|}.$$

As for the numerator in (43), define $g_\varepsilon(x,s) := x^2 + \varepsilon^2 + 2x\varepsilon\cos(s)$ and note that $|\int_{\mathbb{R}^+} [\varepsilon e^{\mathrm{i}s} + x]^{\alpha-1} \eta(dx)| \leq \int_{\mathbb{R}^+} [g_\varepsilon(x,s)]^{(\alpha-1)/2} \leq \varepsilon^{\alpha-1}$ since $\alpha < 1$ and, as before, $|s| < \pi/2$. On the other hand, if we further set $d_\varepsilon(x,s) := \arctan(\varepsilon\sin(s)/[x + \varepsilon\cos(s)])$, then

$$\left| \int_{\mathbb{R}^+} [\varepsilon e^{\mathrm{i}s} + x]^\alpha \eta(dx) \right| \geq \int_{\mathbb{R}^+} [g_\varepsilon(x,s)]^{\alpha/2} \cos(\alpha d_\varepsilon(x,s)) \eta(dx).$$

For any $s \in (-\frac{\pi}{2}, \frac{\pi}{2})$ and $x \geq 0$, it can be seen that $\cos(\alpha d_\varepsilon(x,s)) \geq \cos(\alpha\pi/2)$ and $g_\varepsilon(x,s) \geq x^2$. These, in turn, yield the following lower bound for the denominator in (43):

$$\left| \left\{ \int_{\mathbb{R}^+} [\varepsilon e^{\mathrm{i}s} + x]^\alpha \eta(dx) \right\}^{\theta/\alpha+1} \right| \geq \left\{ \cos\left(\frac{\alpha\pi}{2}\right) \int_{\mathbb{R}^+} x^\alpha \eta(dx) \right\}^{\theta/\alpha+1} = K_{\alpha,\theta} > 0.$$

Hence, the expression in (43) is, for any $\varepsilon > 0$ and $(x,s) \in \mathbb{R}^+ \times (-\frac{\pi}{2}, \frac{\pi}{2})$, bounded by some constant depending only on $\alpha$ and $\theta$. This implies (42) and then

$$\int_{\Gamma_y} h(w)\, dw = \lim_{\varepsilon \downarrow 0} \int_0^y [h(-x - i\varepsilon) - h(-x + i\varepsilon)]\, dx.$$

In order to interchange the limit with the integral, we now show that for any $y$, there exists $B_y \subset (0, y)$, with $\lambda(B_y^c) = 0$, such that for any $x \in B_y$, one has $|h(-x - i\varepsilon)| \leq \bar{h}(x)$ and $\bar{h}$ is integrable on $(0, y)$. The same argument can be applied to $|h(-x + i\varepsilon)|$. If we set

$$(44) \quad \gamma_{\varepsilon,\alpha}(x) := \int_{\mathbb{R}^+} [(t-x)^2 + \varepsilon^2]^{\alpha/2} \cos\left(\alpha \arctan \frac{\varepsilon}{t-x} + \alpha\pi \mathbb{I}_{(0,x)}(t)\right) \eta(dt),$$

$$(45) \quad \zeta_{\varepsilon,\alpha}(x) := \int_{\mathbb{R}^+} [(t-x)^2 + \varepsilon^2]^{\alpha/2} \sin\left(\alpha \arctan \frac{\varepsilon}{t-x} + \alpha\pi \mathbb{I}_{(0,x)}(t)\right) \eta(dt),$$

then

$$|h(-x - i\varepsilon)| = |y - x - i\varepsilon|^{\theta-1} \frac{[\gamma_{\varepsilon,\alpha-1}^2(x) + \zeta_{\varepsilon,\alpha-1}^2(x)]^{1/2}}{[\gamma_{\varepsilon,\alpha}^2(x) + \zeta_{\varepsilon,\alpha}^2(x)]^{(\theta+2\alpha)/(2\alpha)}}.$$

Now, note that $[\gamma_{\varepsilon,\alpha-1}^2(x) + \zeta_{\varepsilon,\alpha-1}^2(x)]^{1/2} \leq K \int_{\mathbb{R}^+} |x-t|^{\alpha-1} \eta(dt)$, for some constant $K > 0$. Furthermore, for any $t \in [x, +\infty)$, one has $\alpha \arctan(\varepsilon/(t-x)) \in (0, \alpha\pi/2]$ and for any $t \in (0, x)$ it follows that $\alpha\pi + \alpha \arctan(\varepsilon/(t-x)) \in (\alpha\pi/2, \alpha\pi)$. Hence, if $\alpha \in (0, 1/2]$, then $\gamma_{\varepsilon,\alpha}^2(x) \geq M > 0$ for any $x$. On the other hand, if $\alpha \in (1/2, 1)$, then $\sin(\alpha \arctan \varepsilon/(t-x)) \leq \sin(\arctan \varepsilon/(t-x)) = \varepsilon/\sqrt{\varepsilon^2 + (t-x)^2}$ for any $t \geq x$ and

$$\zeta_{\varepsilon,\alpha}(x) = \int_{[x,+\infty)} |t - x - i\varepsilon|^\alpha \sin\left(\alpha \arctan \frac{\varepsilon}{t-x}\right) \eta(dt)$$

$$- \int_{(0,x)} |t - x - i\varepsilon|^\alpha \sin\left(\alpha\pi + \alpha \arctan \frac{\varepsilon}{t-x}\right) \eta(dt)$$

$$\leq \varepsilon \int_{[x,+\infty)} [(t-x)^2 + \varepsilon^2]^{(\alpha-1)/2} \eta(dt) - \sin(c_\alpha) \int_{(0,x)} (x-t)^\alpha \eta(dt)$$

$$\leq \varepsilon^\alpha - \sin(c_\alpha) \int_{(0,x)} (x-t)^\alpha \eta(dt),$$

where $c_\alpha = \arg\min_{\alpha\pi/2 < c < \alpha\pi} \sin(c)$. Since $x \in C_\eta$, there exists $\varepsilon^* > 0$ such that, for any $\varepsilon \in (0, \varepsilon^*)$, one has $\zeta_{\varepsilon,\alpha}(x) \leq M_\alpha < 0$. This implies that $|h(-x - i\varepsilon)| \leq M'[y^2 + M'']^{(\theta-1)/2} \int_{\mathbb{R}^+} |x-t|^{\alpha-1} \eta(dt)$ for some suitable positive constants $M'$ and $M''$. The proof is then complete if we can show that $x \mapsto$

$\int_{\mathbb{R}^+} |x-t|^{\alpha-1} \eta(dt)$ is integrable on $(0, y)$. To this end, note that

$$\int_0^y \int_{\mathbb{R}^+ \setminus \{x\}} |x-t|^{\alpha-1} \eta(dt)\, dx$$

$$= \int_{(0,y)} \left\{ \int_0^t (t-x)^{\alpha-1}\, dx + \int_t^y (x-t)^{\alpha-1}\, dx \right\} \eta(dt)$$

$$+ \int_{(y,+\infty)} \int_0^y (t-x)^{\alpha-1}\, dx\, \eta(dt)$$

and this turns out to be finite since $f \in \mathscr{E}_\alpha(H)$ yielding $\int t^\alpha \eta(dt) < \infty$.

## APPENDIX B: PROOFS FOR SECTION 4

We now prove the main results stated in Section 4 concerning the determination of the probability distribution of the mean of a normalized $\alpha$-stable subordinator.

**B.1. Proof of Theorem 4.1.** The first thing to note is that

$$(46) \qquad \mathcal{S}_1[z; M_{\alpha,0}(\eta)] = \frac{\int [z+x]^{\alpha-1} \eta(dx)}{\int [z+x]^\alpha \eta(dx)}$$

for any $z$ such that $|\arg(z)| < \pi$. In order to evaluate the density $q_{\alpha,\eta}$, one can invert (46) by means of the Perron–Stieltjes formula, which yields

$$q_{\alpha,0}(y) = \frac{1}{2\pi i} \lim_{\varepsilon \downarrow 0} \{ \mathcal{S}_1[-y - i\varepsilon; M_{\alpha,0}(\eta)] - \mathcal{S}_1[-y + i\varepsilon; M_{\alpha,0}(\eta)] \}$$

and it can be seen that the above reduces to

$$q_{\alpha,0}(y) = \frac{1}{\pi} \lim_{\varepsilon \downarrow 0} \mathrm{Im}\{ \mathcal{S}_1[-y - i\varepsilon; M_{\alpha,0}(\eta)] \} = \frac{1}{\pi} \lim_{\varepsilon \downarrow 0} \mathrm{Im} \frac{\gamma_{\varepsilon,\alpha-1}(y) - i\zeta_{\varepsilon,\alpha-1}(y)}{\gamma_{\varepsilon,\alpha}(y) - i\zeta_{\varepsilon,\alpha}(y)}.$$

The assumptions $\int_{\mathbb{R}^+} x^\alpha \eta(dx) < \infty$ and $y \in \Theta_{\alpha,\eta}$ allow a straightforward application of the dominated convergence theorem. This leads to $\lim_{\varepsilon \downarrow 0} \gamma_{\varepsilon,\alpha}(y) = \gamma_\alpha(y)$ and $\lim_{\varepsilon \downarrow 0} \zeta_{\varepsilon,\alpha}(y) = \zeta_\alpha(y)$ for any $y$, while $\lim_{\varepsilon \downarrow 0} \gamma_{\varepsilon,\alpha-1}(y) = \gamma_{\alpha-1}(y)$ and $\lim_{\varepsilon \downarrow 0} \zeta_{\varepsilon,\alpha-1}(y) = \zeta_{\alpha-1}(y)$ for any $y \in \Theta_{\alpha,\eta}$. The result then easily follows.

**B.2. Proof of Theorem 4.2.** From Theorem 4.1, it is known that a density function $q_{\alpha,0}$ for $\int_{\mathbb{R}^+} x \tilde{P}_{\alpha,0}(dx)$ is of the form

$$q_{\alpha,0}(y) = \Delta_{\alpha,1}(y) = \frac{\gamma_{\alpha-1}(y)\zeta_\alpha(y) - \gamma_\alpha(y)\zeta_{\alpha-1}(y)}{\pi\{\gamma_\alpha^2(y) + \zeta_\alpha^2(y)\}}$$

for any $y \in \Theta_{\alpha,\eta}$. Let us now compute the derivatives of $\gamma_\alpha$ and $\zeta_\alpha$. In order to do so, note that

$$\lim_{h \downarrow 0} \frac{\int_{(0,y+h)}(y+h-x)^\alpha \eta(dx) - \int_{(0,y)}(y-x)^\alpha \eta(dx)}{h}$$
$$= \lim_{h \downarrow 0} \int_{(0,y)} \frac{(y+h-x)^\alpha - (y-x)^\alpha}{h} \eta(dx)$$
$$- \lim_{h \downarrow 0} \int_{[y,y+h)} \frac{(y+h-x)^\alpha}{h} \eta(dx).$$

Since we are confining ourselves to evaluating the density on the set of points $y$ in $\Theta_{\alpha,\eta}$, the probability measure $\eta$ does not have a positive mass on such $y$'s. Hence, for suitable positive constants $c_y$ and $k_y$, one has

$$\int_{[y,y+h)} \frac{(y+h-x)^\alpha}{h} \eta(dx) = c_y h^{\alpha-1} \{\Psi(y+h) - \Psi(y)\} \leq k_y h^\alpha,$$

where the first equality follows from the mean value theorem for Riemann–Stieltjes integrals and the inequality is a consequence of the fact that $\Psi$ is Lipschitz of order 1 at $y$. On the other hand, $\{(y+h-x)^\alpha - (y-x)^\alpha\}/h \leq \alpha(y-x)^{\alpha-1}$ for any $x \in (0,y)$ and for any $h > 0$. Since $x \mapsto (y-x)^{\alpha-1}$ is integrable on $(0,y)$ for any $y \in \Theta_{\alpha,\eta}$, one can apply the dominated convergence theorem to obtain $(d/dy) \int_{(0,y)} (y-x)^\alpha \eta(dx) = \alpha \int_{(0,y)} (y-x)^{\alpha-1} \eta(dx)$. The same argument can be applied to prove that $(d/dy) \int_{(y,+\infty)} (x-y)^\alpha \eta(dx) = -\alpha \int_{(y,+\infty)} (x-y)^{\alpha-1} \eta(dx)$. These imply that $\zeta_{\alpha-1}(y) = -\alpha^{-1} \zeta'_\alpha(y)$ and $\gamma_{\alpha-1}(y) = -\alpha^{-1} \gamma'_\alpha(y)$. Consequently,

$$q_{\alpha,0}(y) = \frac{1}{\alpha \pi} \frac{d}{dy} \arctan \frac{\zeta_\alpha(y)}{\gamma_\alpha(y)},$$

from which one easily obtains the expression of the c.d.f. $x \mapsto \int_0^x q_{\alpha,0}(y)\, dy$ displayed in the statement of the theorem.

## APPENDIX C: PROOFS FOR SECTION 5

**C.1. Proof of Theorem 5.1.** From the definition of $\Delta_{\alpha,\theta}$ and the representation of the generalized Stieltjes transform of $M_{\alpha,\theta}(\eta)$, as given in [44] and [45], it is apparent that

$$\Delta_{\alpha,\theta}(t) = \frac{1}{\pi} \lim_{\varepsilon \downarrow 0} \operatorname{Im} \mathcal{S}_\theta[-t - i\varepsilon; M_{\alpha,\theta}(\eta)]$$
$$= \frac{1}{\pi} \lim_{\varepsilon \downarrow 0} \operatorname{Im} \left\{ \int_\mathbb{X} (-t - i\varepsilon + x)^\alpha \eta(dx) \right\}^{-\theta/\alpha},$$

where $\mathbb{X} \subset \mathbb{R}^+$. One has

$$\left\{\int_{\mathbb{R}} (-t - i\varepsilon + x)^\alpha \eta(dx)\right\}^{-\theta/\alpha} = \exp\left\{-\frac{\theta}{\alpha} \log(\gamma_{\varepsilon,\alpha}(t) - i\zeta_{\varepsilon,\alpha}(t))\right\}.$$

Let us first confine our attention to the case in which $\alpha$ is in the interval $(0, 1/2]$. Since $\alpha \arctan(\frac{\varepsilon}{x-t}) + \alpha\pi \mathbb{I}_{(0,t)}(x) \in (0, \alpha\pi)$, for any $t$ and $x$, one has $\zeta_{\varepsilon,\alpha}(t) > 0$ and $\gamma_{\varepsilon,\alpha}(t) > 0$. Consequently,

$$\exp\left\{-\frac{\theta}{\alpha} \log(\gamma_{\varepsilon,\alpha}(t) - i\zeta_{\varepsilon,\alpha}(t))\right\}$$
$$= \{\gamma_{\varepsilon,\alpha}^2(t) + \zeta_{\varepsilon,\alpha}^2(t)\}^{-\theta/(2\alpha)} \exp\left\{i\frac{\theta}{\alpha} \arctan \frac{\zeta_{\varepsilon,\alpha}(t)}{\gamma_{\varepsilon,\alpha}(t)}\right\}.$$

Note that the absolute values of each of the two integrands defining $\gamma_{\varepsilon,\alpha}$ and $\zeta_{\varepsilon,\alpha}$ are bounded by $|x - t|^\alpha + K$, which is integrable with respect to $\eta$. We can then apply a dominated convergence argument to obtain

$$\lim_{\varepsilon \downarrow 0} \gamma_{\varepsilon,\alpha}(t) = \gamma_\alpha(t), \qquad \lim_{\varepsilon \downarrow 0} \zeta_{\varepsilon,\alpha}(t) = \zeta_\alpha(t)$$

for any $t > 0$. This implies (30) after noting that, in this case, $\Gamma_\alpha = \varnothing$.

On the other hand, when $\alpha \in (1/2, 1)$, one needs to consider the set $\Gamma_{\varepsilon,\alpha} := \{t \in \mathbb{R}^+ : \gamma_{\varepsilon,\alpha}(t) > 0\}$ and note that $\Gamma_{\varepsilon,\alpha}^c \cap (0, y)$ is nonempty for some values of $y$ in $C_\eta$. This yields a slightly different form for the arguments of the complex numbers involved in the definition of $\Delta_{\alpha,\theta}$. One can easily mimic the line of reasoning employed for the case $\alpha \in (0, 1/2]$ so as to again obtain (30).

**C.2. Proof of Theorem 5.2.** Since $M_{\alpha,\theta}(\eta) \stackrel{d}{=} M_{0,\theta}(Q_{\alpha,0})$, it follows that the distribution functions given in (10) and (24) are equal for all $\theta > 0$. Statement (i) then follows by the uniqueness properties of the integral representations. The first identity in statement (ii) follows immediately by setting $\theta = \alpha$ in statement (i) which, noting that $0 < \alpha < 1$, uses the strict positivity $\sin(\pi\alpha Q_{\alpha,0}(t))$ for $Q_{\alpha,0}(t) > 0$. If one now exploits (i) with $\theta = \alpha$, it is possible to obtain

$$e^{-R_\alpha(t)} = \left[\frac{\Delta_{\alpha,\alpha}(t)\pi}{\sin(\pi\alpha Q_{\alpha,0}(t))}\right]^{1/\alpha}.$$

We now set $\theta = \alpha$ in (32) to obtain the second identity involved in (ii), thus completing the proof.

**C.3. Proof of Theorem 5.3.** By definition,

$$\Delta_{\alpha,\theta+1}(t) = \frac{1}{\pi} \lim_{\varepsilon \downarrow 0} \{\mathcal{S}_{\theta+1}[-t - i\varepsilon; M_{\alpha,\theta}(\eta)] - \mathcal{S}_{\theta+1}[-t + i\varepsilon; M_{\alpha,\theta}(\eta)]\},$$

which can be seen to imply

$$\Delta_{\alpha,\theta+1}(t) = \frac{1}{\pi} \lim_{\varepsilon \downarrow 0} \operatorname{Im} \frac{\int_{\mathbb{R}} (-t - i\varepsilon + x)^{\alpha-1} \eta(dx)}{\{\int_{\mathbb{R}} (-t - i\varepsilon + x)^{\alpha} \eta(dx)\}^{(\theta+\alpha)/\alpha}}.$$

For any $\varepsilon > 0$, $|(-t - i\varepsilon + x)^\alpha|$ can be bounded by an integrable function with respect to $\eta$, not depending on $\varepsilon$ in a similar fashion as in the proof of Theorem 4.1. On the other hand, $|(-t - i\varepsilon + x)^{\alpha-1}| \leq |x - t|^{\alpha-1} + K'$, for any $\varepsilon > 0$, $x$ and $t$. If we further set $t \in \Theta_{\alpha,\eta}$, then $x \mapsto |x - t|^{\alpha-1}$ is integrable with respect to $\eta$ and the dominated convergence theorem can be applied. The expression in (33) easily follows.

**C.4. Proof of Proposition 5.1.** Statement (i) concerns the evaluation of $\mathbb{E}_\eta[\omega_m^{m\kappa}(t) \sin(\pi m \kappa \eta_m(t))]$. Here, we use the fact that if $c_{m\kappa}(t) < \infty$, there exists, by a change of measure, a density for each $Y_k$ which is proportional to $|t - y|^{-m\kappa} \eta(dy)$. It then follows that, with respect to this i.i.d. law for $(Y_k)$, $m\eta_m(t)$ is a Binomial($mp_{m\kappa}(t)$) random variable and the result is proved. Statement (ii) is derived from (8) using a conditioning argument. Similarly, statement (iii) follows from (10), noting that the jumps of $\theta \eta_m$ are less than 1.

## APPENDIX D: PROOF OF THEOREM 6.1

The proof follows as a direct consequence of the mixture representation of the laws of the $\tilde{P}_{\alpha,\theta}$ random probability measures deduced from their posterior distributions. Specifically, one can immediately deduce from [33], with $n = 1$, that,

$$\tilde{P}_{\alpha,\theta}(\cdot) \stackrel{d}{=} B_{\theta+\alpha,1-\alpha} \tilde{P}_{\alpha,\theta+\alpha}(\cdot) + (1 - B_{\theta+\alpha,1-\alpha}) \delta_Y(\cdot),$$

yielding the stated result. Specifically, apply the above identity to $\tilde{P}_{\alpha,\theta}(g)$, where $g(x) = x$. Naturally, this statement is an extension of the result deduced from Ferguson's characterization of a posterior distribution of a Dirichlet process; see [10, 11], also related discussions about mixture representations derived from posterior distributions in [21, 22]. As for the proof of (ii), this follows from the distributional identity stated in Theorem 2.1 and by setting $Y = M_{\alpha,0}(\eta)$ in (i). Indeed, $M_{\alpha,\theta}(\eta) \stackrel{d}{=} M_{0,\theta}(Q_{\alpha,0}) \stackrel{d}{=} B_{\theta,1} M_{0,\theta}(Q_{\alpha,0}) + (1 - B_{\theta,1}) M_{\alpha,0}(\eta)$. The proof is complete by again applying Theorem 2.1 to the first summand of the last sum in the previous chain of identities.

**Acknowledgments.** The authors are grateful to an Associate Editor and an anonymous referee for their valuable comments and suggestions.

L. F. JAMES  
DEPARTMENT OF INFORMATION  
AND SYSTEMS MANAGEMENT  
HONG KONG UNIVERSITY OF  
SCIENCE AND TECHNOLOGY  
CLEAR WATER BAY, KOWLOON  
HONG KONG  
E-MAIL: lancelot@ust.hk

A. LIJOI  
DIPARTIMENTO DI ECONOMIA POLITICA  
E METODI QUANTITATIVI  
UNIVERSITÀ DEGLI STUDI DI PAVIA  
VIA SAN FELICE 5  
27100 PAVIA  
ITALY  
E-MAIL: lijoi@unipv.it

I. PRÜNSTER  
DIPARTIMENTO DI STATISTICA  
E MATEMATICA APPLICATA  
UNIVERSITÀ DEGLI STUDI DI TORINO  
PIAZZA ARBARELLO 8  
10122 TORINO  
ITALY  
E-MAIL: igor@econ.unito.it